\documentclass[10pt,a4paper]{article}
\usepackage{amssymb, amsmath,amsfonts,amsthm}
\usepackage{graphicx}
\usepackage{color}
\usepackage[all]{xy}
\usepackage{psfrag}

\theoremstyle{plain}
\newtheorem{theorem}{Theorem}[section]
\newtheorem{proposition}[theorem]{Proposition}
\newtheorem{corollary}[theorem]{Corollary}
\newtheorem{claim}{Claim}
\newtheorem{lemma}[theorem]{Lemma}

\theoremstyle{definition}

\newtheorem{example}[theorem]{Example}

\newcommand{\R}{\mathbb{R}}

\newcommand{\Z}{\mathbb{Z}}
\newcommand{\N}{\mathbb{N}}
\newcommand{\UT}{\mathrm{UT}}

\DeclareMathOperator{\Aut}{Aut} \DeclareMathOperator{\Aff}{Aff}
\DeclareMathOperator{\ATF}{ATF} \DeclareMathOperator{\ITF}{ITF}
 
\DeclareMathOperator{\SL}{SL} \DeclareMathOperator{\C}{\mathbb{C}}
 \DeclareMathOperator{\Rad}{Rad}
\DeclareMathOperator{\Diam}{Diam} 
\DeclareMathOperator{\Hom}{Hom} \DeclareMathOperator{\id}{id}
\DeclareMathOperator{\Out}{Out}
\DeclareMathOperator{\Hol}{Hol}
\DeclareMathOperator{\BS}{BS}  
\newcommand{\pf}{\textbf{Proof: }}
\renewcommand{\qed}{\hfill{$\square$}}

\newcommand{\be}{\begin{enumerate}}
\newcommand{\ee}{\end{enumerate}}

\renewcommand{\a}{\mathbf{a}}
\renewcommand{\b}{\mathbf{b}}\renewcommand{\c}{\mathbf{c}}
\newcommand{\h}{\hspace*{0.4cm}}
\newcommand{\Y}{\mathrm{Y}}
\newcommand{\Br}{\mathrm{Br}}

\begin{document}

\markboth{SHANE O ROURKE} {Affine actions on non-archimedean trees}

\bibliographystyle{plain}

\author{S. O Rourke}
\title{Affine actions on non-archimedean trees}
\maketitle

\begin{abstract}
We initiate the study of affine actions of groups on $\Lambda$-trees
for a general ordered abelian group\ $\Lambda$; these are actions by
dilations rather than isometries. This gives a common generalisation
of isometric action on a $\Lambda$-tree, and affine action on an
$\R$-tree as studied by I. Liousse. The duality between based length
functions and actions on $\Lambda$-trees is generalised to this
setting. We are led to consider a new class of groups: those that
admit a free affine action on a $\Lambda$-tree for some $\Lambda$.
Examples of such groups are presented, including soluble
Baumslag-Solitar groups and the discrete Heisenberg group.
\end{abstract}

\section*{Introduction}

This paper attempts to draw together three strands: the study of
affine actions on $\R$-trees, which were first considered by I.
Liousse in \cite{Liousse}; Lyndon length functions (based length
functions) on groups and the duality between such length functions
and isometric actions on $\Lambda$-trees; and a desire to
rehabilitate certain examples of groups (studied in \cite{thesis})
that admit a free non-nesting action on a pretree but which admit no
free isometric action on a $\Lambda$-tree.

Isometric actions on $\Lambda$-trees now form a well-established
theme in geometric group theory. This area has its origins in a
paper of R. Lyndon \cite{Lyndon}, where Nielsen's proof of the
classical Nielsen-Schreier Theorem was presented in terms of
abstract length functions on groups, and the Kurosh Subgroup Theorem
was also deduced from this framework. In fact, anticipating future
developments in the area, Lyndon commented that ``we feel strongly
that the restriction to length functions whose values are integers,
rather than real numbers, elements of an ordered abelian group, or
even of algebraic structures of a more unfamiliar nature is
regrettable\ldots''

A crucial step in the development of the theory was the development
of Bass-Serre theory in the late 1960s. This theory established an
equivalence between group actions (without inversions) on trees and
graphs of groups (the latter are a common generalisation of
amalgamated free products and HNN extensions), and provided new
proofs of many results and a new viewpoint enabling the development
of others.

The link between group actions on trees and length functions was
provided by I.M. Chiswell \cite{Chiswell-length}, who showed that
specifying a Lyndon length function on a group is essentially
equivalent to giving an action of the group on a tree. One now had
three strikingly different conceptual approaches to a common
phenomenon: length functions, actions on trees, and graphs of
groups. This framework had the effect of bringing many geometric and
group theoretic themes closer together. For example these ideas gave
a common conceptual viewpoint for understanding the two original
proofs of the Nielsen-Schreier subgroup theorem.

The theory of $\R$-trees and their interaction with group theory has
become a very important theme over the last three decades. We refer
to \cite{Bestvina-survey} for more information about this rich
branch of the subject. Suffice it to say that $\Lambda$-trees, and
with them  $\R$-trees, were first defined by J. Morgan and P. Shalen
in the course of their work on algebraic subvarieties of
$\SL(2,\C)$-characters of finitely generated groups. (It is worth
mentioning however that the notion of $\R$-trees had been implicit
beforehand in the work of W. Imrich \cite{Imrich} and J. Tits
\cite{Tits}, as well as the work of Chiswell cited above.)

Free actions have always been of special interest. A group that
admits a free action (without inversions) on a $\Lambda$-tree is
said to be $\Lambda$-free. This terminology stems from the fact that
in the classical case $\Lambda=\Z$, a group is $\Lambda$-free if and
only if it is free. The question of which groups are $\R$-free was
the subject of intense investigation culminating in Rips' Theorem: a
finitely generated group is $\R$-free if and only if it is a free
product of finitely generated free abelian groups and surface groups
(with the exception of the fundamental groups of non-orientable
surfaces of genus $\leq 3$).

It is of interest to note here that finitely generated fully
residually free groups are $\Lambda$-free, where $\Lambda$ can be
taken to have the form $\Z^n$. Moreover $\R$-free groups are locally
fully residually free, as a consequence of Rips' Theorem. However
$\Z^n$-free groups are not typically residually free: the class of
$\Z^n$-free groups is thus considerably larger than the class of
fully residually free groups.

H. Bass \cite{Bass} examined group actions on
$\Z\times\Lambda_0$-trees, which in principle makes group actions on
$\Z^n$-trees tractable. Recently, O. Kharlampovich, A. Myasnikov, V.
Remeslennikov and D. Serbin have given a detailed description of
$\Z^n$-free groups \cite{KM-Zn} (see also \cite{KM-Lambda} and
\cite{KM-inf-words}). Moreover they have shown that in the finitely
presented case, if a group $G$ admits a regular free action on a
$\Lambda$-tree for some ordered abelian group\ $\Lambda$, then $G$
admits a regular free action on an $\R^n$-tree for some $n$.

A new ingredient to the theory was introduced into the theory by
Myasnikov, Remeslennikov and Serbin, namely infinite $\Lambda$-words
in a group. They showed that a group that admits a faithful
representation by infinite $\Lambda$-words also admits a free action
on a $\Lambda$-tree \cite{inf-words-free-actions}. Conversely, a
result of Chiswell \cite{free-actions-inf-words} showed that a group
equipped with a free action on a $\Lambda$-tree admits a faithful
representation as infinite $\Lambda$-words.

Of course the appearance of generalised words in the theory should
come as no surprise since the prototypical example of a free action
on a ($\Z$-)tree is the action of a free group on its Cayley graph
(with respect to a free generating set). The Lyndon lengths of
elements of the group with respect to this action coincide with the
lengths as words in the free generating set. A similar statement is
true of the general case.

Further progress has been made by A. Nikolaev and Serbin
\cite{Nik-Serb} who have given an effective solution of the
membership problem and the power problem for groups that act freely
on $\Z^n$-trees.

The surge in interest in actions on $\Lambda$-trees in recent times
follows the work of Kharlampovich and Myasnikov on the Tarski
problems concerning the elementary theory of non-abelian free
groups. They showed among many other things that any pair of
non-abelian free groups are elementarily equivalent, and that the
elementary theory of such groups is decidable. These are also called
limit groups, a term introduced by Z. Sela \cite{Sela-1} who in a
series of papers gave another solution of Tarski's problems using a
different approach.

A. Martino and the author surveyed the topic of free actions on
$\Z^n$-trees in \cite{Zn-survey} though this predates many of the
developments in the theory outlined above.

A further advance was obtained by V. Guirardel \cite{Guirardel} who
investigated group actions on $\R^n$-trees, and showed that
$\R^n$-free groups admit a graph of groups decomposition with
$\R^{n-1}$-free vertex groups. He also showed that such groups are
coherent, thereby giving a proof of the finite presentability of
limit groups.\\

All actions referred to so far are isometric actions. Continuous
actions on $\R$-trees by non-nesting automorphisms have also been
considered by G. Levitt \cite{Levitt}, and the more general
situation of non-nesting actions on pretrees was studied by B.
Bowditch and J. Crisp \cite{Bowditch-Crisp}. In particular the
latter have shown that a finitely presented group that admits a
non-trivial archimedean action on a median pretree also admits a
continuous non-trivial non-nesting action on an $\R$-tree. Thus by
the main result of \cite{Levitt}, such a group admits a non-trivial
isometric action on an $\R$-tree.

In \cite{Liousse}, Liousse initiated the study of affine actions on
$\R$-trees. An affine automorphism of an ($\R$-)metric space is a
surjective function $\phi:X\to X$ for which there exists a constant
$\alpha_\phi$ such that $d(\phi x,\phi y)=\alpha_\phi d(x,y)$ for
all $x,y\in X$. (Of course such an $\alpha_\phi$ must be positive
for non-degenerate $X$.) Examples of affine actions include the
following, due to F. Paulin \cite{Paulin-affine}. Suppose that
$\Gamma$ is a hyperbolic group, and $H$ a subgroup of $\Aut(\Gamma)$
with amenable image in $\Out(\Gamma)$  and with infinite centre.
Then $\Gamma\rtimes H$ admits an affine action on an $\R$-tree such
that the restriction to $\Gamma$ is isometric, has no global fixed
point, and has virtually cyclic arc stabilisers.

Liousse gave two families of examples of groups (recalled as
Example~\ref{liousse-exx} below) that admit free affine actions on
$\R$-trees but which admit no free isometric action on an $\R$-tree.
Martino and the author \cite{MOR} later showed that Liousse's groups
do admit free isometric actions on $\Z^n$-trees for some $n$. This
work may be thought of as an attempt to admit affine actions on
$\R$-trees to the fold of isometric actions on $\Lambda$-trees at
the expense of modifying the ordered abelian group. However, it is
also natural to reverse this process --- that is, to attempt to
broaden the scope of affine actions on $\R$-trees to incorporate
non-archimedean $\Lambda$.

In \cite{thesis} the author considered the problem of equivariantly
embedding a pretree in a $\Lambda$-tree equipped with an isometric
action, in the presence of a group action on the former. One
necessary condition for such an embedding to exist is
\emph{rigidity}: a map $g$ is rigid if $g[x,y]\not\subset[x,y]$ and
$g[x,y]\not\supset[x,y]$ for all points $x$ and $y$. However it was
shown that this condition is not sufficient as the wreath product
$C_\infty\wr C_\infty$ of two infinite cyclic groups admits a free,
rigid action on a pretree but does not admit a free isometric action
on any $\Lambda$-tree. From the point of view of isometric actions,
$C_\infty\wr C_\infty$ may thus appear to be somewhat pathological.

However, the concept of affine action can be naturally extended to
$\Lambda$-trees. The key idea here is to require $d(gx,gy)=\alpha_g
d(x,y)$ where now $\alpha_g$ is an element of $\Aut^+(\Lambda)$, the
group of order-preserving group automorphisms of $\Lambda$.
In this framework, isometric actions on $\Lambda$-trees and affine
actions on $\R$-trees in the sense of Liousse appear as special
cases. Moreover, examples such as the action of $C_\infty\wr
C_\infty$ mentioned above --- and many others besides --- fall under
this new heading of affine action on a $\Lambda$-tree. The examples
presented here may be summarised as follows.

\begin{theorem}\label{main-exx}
\be
\item The Heisenberg group $\UT(3,\Z)$ admits a free affine action on a $\Z^3$-tree.
\item The wreath product of two infinite cyclic groups admits a free affine action on a $\Z\times\R$-tree.
\item The soluble Baumslag-Solitar groups $\BS(1,a)$ admit a free affine action on a $\Z\times\R$-tree.
\ee
\end{theorem}

The main goal of the current paper is to commence a systematic study
of affine actions on $\Lambda$-trees $X$, where $\Lambda$ is not
necessarily archimedean --- the titular non-archimedean trees. We
will focus especially on constructions that yield free affine
actions. In future work, we propose to focus on the case
$\Lambda=\Z^n$, and on the case where $X$ is a linear tree. Various
technical difficulties present themselves, notably the absence of a
good analogue of hyperbolic length function. Nevertheless much of
our work follows the general programme of the theory of isometric
actions. For example, we show that the duality between Lyndon length
functions and isometric actions on $\Lambda$-trees, which was first
proven by Chiswell \cite{Chiswell-length} in the case $\Lambda=\Z$
extends naturally to the affine case.

\begin{theorem}\label{LLF-duality}
Let $\alpha:G\to\Aut^+(\Lambda)$ be a homomorphism, fix an
$\alpha$-affine action of $G$ on a $\Lambda$-tree $X$, and a point
$x\in X$. The function defined by $L(g)=L_x(g)=d(x,gx)$ is an
\emph{$\alpha$-affine Lyndon length function} on $G$.

Conversely, for any $\alpha$-affine Lyndon length function $L$ on
$G$ there is an (essentially unique) $\alpha$-affine action of $G$
on a $\Lambda$-tree, and a basepoint $x$ such that $L=L_x$.
\end{theorem}

We also generalise some of Liousse's results, for example showing
that under certain conditions an isometric action of a normal
subgroup of $G$ can be extended to an affine action of $G$. Some of
the results established in the general case give rise to results
that are new even in the archimedean case. For example, we show

\begin{theorem}\label{main-free-product}
The class of groups that admit a free affine action on a
$\Lambda$-tree is closed under free products.\end{theorem}

One of the striking features of the theory of free affine actions on
non-archimedean trees is that stabilisers of lines are no longer
necessarily torsion-free abelian. In the case $\Lambda=\Z^n$, a
group that stabilises a line and preserves the orientation is
finitely generated torsion-free nilpotent, while for $\Lambda$ with
$\Aut^+(\Lambda)$ soluble, the line stabilisers are torsion-free
soluble, but not necessarily nilpotent. This means that certain
familiar properties of isometrically `tree-free' groups such as the
CSA property must be generalised.

Recently, we have shown that groups that admit free affine actions
on linear $\Lambda$-trees form a much larger class than would be
predicted from the isometric case. Namely, \emph{residually
torsion-free nilpotent groups} admit free affine actions on linear
$\Lambda$-trees. This is shown in \cite{aff-linear-case}. It is, as
far as we are aware, an open question whether groups that admit a
free isometric action on a $\Z^n$-tree are residually torsion-free
nilpotent.

I would like to thank Martin Bridson for some helpful comments
concerning the example of Theorem~\ref{main-exx}(1). I would also
like to thank Enric Ventura and the other organisers of the fifth
conference on Geometric and Asymptotic Group Theory with
Applications held in Manresa in July 2011 for such a stimulating and
enjoyable conference.

\section{Fundamentals}

\subsection{Ordered abelian groups and their automorphism groups}

Let $\Lambda$ (or more precisely $(\Lambda,+,\leq)$) be an ordered
abelian group. We write $\Aut^+(\Lambda)$ for the group of all group
automorphisms of $\Lambda$ that preserve the order. A \emph{convex}
subgroup of $\Lambda$ is a subgroup $\Lambda_0$, for which
$a,c\in\Lambda_0$ with $a\leq b\leq c$ implies $b\in\Lambda_0$. The
set of convex subgroups of an ordered abelian group\ $\Lambda$ is
linearly ordered by inclusion, and any $\alpha_*\in\Aut^+(\Lambda)$
induces an order-preserving permutation of the convex subgroups of
$\Lambda$. We call an automorphism \emph{shift-free} if this induced
permutation is trivial ---  that is, if it stabilises each convex
subgroup. We write $\Aut_0^+(\Lambda)$ for the group of shift-free
automorphisms of $\Lambda$. We will call an ordered abelian group\
\emph{convex-rigid} if all order-preserving automorphisms are
shift-free.

The order type of the set of non-zero convex subgroups of $\Lambda$
is the \emph{rank} of $\Lambda$. If $\Lambda$ has finite rank $n$
then $\Lambda$ embeds (as an ordered abelian group) in $\R^n$ ---
here $\Lambda_1\times\Lambda_2$ is ordered lexicographically via
$(\lambda_1,\lambda_2)\leq (\lambda_1',\lambda_2')$ if
$\lambda_1<\lambda_1'$ or if $\lambda_1=\lambda_1'$ and
$\lambda_2\leq\lambda_2'$. Inductively this gives an order on
$\Lambda_1\times\cdots\times\Lambda_n$ and, as a special case, on
$\Lambda^n_0$, where the $\Lambda_i$ are given ordered abelian
groups. Whenever an ordered abelian group\ is given as a direct
decomposition or as a (direct) power, it is always the lexicographic
order that we will consider unless the contrary is specifically
stated. It is worth noting however that other orders are possible.
For example, $\Z\times\Z$ may be lexicographically ordered, making
it a rank 2 ordered abelian group, or it may have rank 1, via an
embedding in $\R$. Subgroups of $\R$ (with the natural order) have
rank $\leq 1$ and are characterised by the familiar archimedean
property.

The structure of $\Aut^+(\Lambda)$ was considered in a paper by P.
Conrad. For our purposes, it suffices to record the following (see
\cite[\S1]{Conrad}).

\begin{proposition}\label{aut12}
If $\Lambda=\Lambda_1\times\cdots\times\Lambda_n$ (with the lexicographic order), then $\Aut_0^+(\Lambda)$ is isomorphic to the group of matrices of the form $$\left(\begin{array}{cccc}\alpha_{n} & h_{n(n-1)} & \cdots & h_{n1}\\
0 & \alpha_{n-1} & \cdots & h_{(n-1)1}\\ \vdots & \vdots & \ddots &
\vdots\\ 0 & 0 & \cdots & \alpha_1
\end{array}\right)$$ where $\alpha_i\in\Aut_0^+(\Lambda_i)$ and $h_{ij}\in\Hom(\Lambda_j,\Lambda_i)$.

In particular $\Aut_0^+(\Lambda_1\times\Lambda_2)\cong
\left(\Aut_0^+(\Lambda_1)\times\Aut_0^+(\Lambda_2)\right)\ltimes\Hom(\Lambda_1,\Lambda_2)$.
(Here the action of $\Aut_0^+(\Lambda_1)\times\Aut_0^+(\Lambda_2)$
is given by  $(\alpha_1,\alpha_2)\cdot h=\alpha_2 h\alpha_1^{-1}$
for $\alpha_i\in\Aut_0^+(\Lambda_i)$ and
$h\in\Hom(\Lambda_1,\Lambda_2)$.)
\end{proposition}

One can deduce from this result that the group $\Aut^+(\Z^n)$ is
isomorphic to the group $\UT(n,\Z)$ of upper triangular matrices
with all diagonal entries equal to 1. Note also that the group
$\mathrm{T}(n,\R)$ of upper triangular matrices with positive units
on the diagonal, embeds in $\Aut^+(\R^n)$ but is not isomorphic to
it on account of the profusion of (non-order-preserving)
homomorphisms $\R\to\R$.

Next we consider $\Aut^+(\Lambda)$ where $\Lambda\leq\R$. In this
case $\alpha_*\in\Aut^+(\Lambda)=\Aut_0^+(\Lambda)$ is determined by
$\alpha_*\cdot 1$, and has the effect of multiplication by the
positive real number $\alpha_*\cdot 1$. One concludes that
$\Aut_0^+(\Lambda)$ embeds in the multiplicative group of positive
real numbers which is isomorphic to $\R=(\R,+)$.

\begin{corollary}\label{sol-f-rank}
Let $\Lambda=\Lambda_1\times\Lambda_2\times\cdots\times\Lambda_n$,
suppose that $\Lambda_i$ is convex-rigid and that
$\Aut^+(\Lambda_i)$ is soluble for all $i$. Then $\Aut^+(\Lambda)$
is soluble. In particular, $\Aut^+(\Lambda)$ is soluble for any
ordered abelian group\ $\Lambda$ of finite rank.
\end{corollary}

Note that in general a group $G$ admits an embedding in
$\Aut^+(\Lambda)$ for some $\Lambda$ precisely when $G$ is right
orderable, a result due to D. Smirnov (see \cite[Theorem
7.1.3]{Mura-Rhemtulla}). In fact R. G\"obel and S. Shelah
\cite{Goebel-Shelah} have shown that $G$ is right orderable if and
only if $G$ is isomorphic to $\Aut^+(\Lambda)$ for some $\Lambda$.

\subsection{Affine length functions and actions on $\Lambda$-trees}

We refer the reader to Chiswell's book \cite{Chiswell-book} where
the basic theory of $\Lambda$-trees is treated in detail. (We state
a characterisation of $\Lambda$-trees just after
Lemma~\ref{length-props} below.)

Let $(X,d)$ and $(X',d')$ be $\Lambda$-trees. In fact most of the
definitions that follow make sense for $\Lambda$-metric spaces in
general. A function $\phi:X\to X'$ is an \emph{$\alpha_*$-affine
map} (or simply an \emph{affine map}) if there exists
$\alpha_*=\alpha_\phi\in\Aut^+(\Lambda)$ such that
$$d'(\phi x,\phi y)=\alpha_*(d(x,y))$$
for all $x,y\in X$. The automorphism $\alpha_*$ is called the
\emph{dilation factor} corresponding to $\phi$. Note that one can
define an affine map somewhat more permissively by requiring only
that $\alpha_*\in\Aut^+(\Lambda_0)$ where $\Lambda_0$ is the
(convex) subgroup generated by the image of the distance function
$d$, but we will not do so here. Note however that if $\Lambda$ has
finite rank then $\alpha_*\in\Aut^+(\Lambda_0)$ extends in a natural
(but typically non-unique) way to
$\bar{\alpha}_*\in\Aut^+(\Lambda)$. Thus $\Aut^+(\Lambda_0)$ embeds
in $\Aut^+(\Lambda)$ in this case.

An \emph{affine isomorphism} is a bijective affine map, and an
\emph{affine automorphism} is an affine isomorphism with the same
domain and codomain. Note that if $\Lambda=\R$ this usage is
consistent with Liousse \cite{Liousse}. An \emph{affine action} of
$G$ on a $\Lambda$-tree $X$ consists of a homomorphism
$\alpha:G\to\Aut^+(\Lambda)$, written $g\mapsto\alpha_g$, together
with an action of $G$ on $X$ such that $g$ is an $\alpha_g$-affine
automorphism of $X$ for all $g\in G$. Such an action may be referred
to as an $\alpha$-affine action. Note that
$\alpha_g$
restricts to the identity map on the image of $d$ precisely when $g$
is an isometry so that the subgroup of $G$ consisting of isometries
contains $\ker\alpha$.

Let $G$ be a group that has an $\alpha$-affine action on a
$\Lambda$-tree $X$, and fix $x_0\in X$. Define a function
$L:G\to\Lambda$ by the formula
$$L(g)=L_{x_0}(g)=d(x_0,gx_0)$$

Such an $L$ is called the \emph{($\alpha$-)based length function}
arising from the action (with respect to $x_0$).

Fix a homomorphism $\alpha:G\to\Aut^+(\Lambda)$. Let $L:G\to\Lambda$
be a function, and for $g,h\in G$, we define the \emph{ancillary
functions} $\a$, $\b$ and $\c$ as follows.

\begin{eqnarray*}\mathbf{a}(g)&=&
\frac{1}{2}\left(L(g)+L(g^{-1})-\alpha_{g^{-1}}L(g^2)\right)\\ \\
\mathbf{b}(g)&=&
\frac{1}{2}\left((\alpha_{g^{-1}}+1)L(g^2)-(\alpha_{g^2}+1)L(g^{-1})\right)\\ \\
\mathbf{c}(g,h)&=&\frac{1}{2}\left(L(g)+L(h)-\alpha_gL(g^{-1}h)\right)\end{eqnarray*}

(Here $1+\alpha_\gamma$ is of course the map $\lambda\mapsto
\lambda+\alpha_\gamma(\lambda)$.)

Call $L$ an \emph{$\alpha$-Lyndon length function} on $G$ if for
$g$, $h$ and $k$ in $G$ one has

\begin{enumerate}
\item[(L1)] $\mathbf{c}(g,h)\in\Lambda$
\item[(L2)] $L(1)=0$
\item[(L3)] $L(g)=\alpha_g (L(g^{-1}))$
\item[(L4)] $(\mathbf{c}(g,h),\mathbf{c}(g,k),\mathbf{c}(h,k))$ is an isosceles triple --- that is, $\c(g,h)<\c(g,k)\Rightarrow\c(g,h)=\c(h,k)$.
\end{enumerate}

If $\alpha_g=1$ for all $g$ then an $\alpha$-Lyndon length function
is a Lyndon length function in the usual sense (apart from (L1),
which is an extra assumption in, for example,
\cite[\S2.4]{Chiswell-book}).

If $L$ is a based length function arising from an $\alpha$-affine
action, then applying the definition of $\mathbf{c}$ to $L$, one
easily sees that $\mathbf{c}(g,h)$ is equal to the length of the
segment $[x,v]=[x,gx]\cap[x,hx]$, as in the isometric case. Property
(L4) of $\alpha$-Lyndon length functions follows, while properties
(L1)-(L3) are immediate. A based length function arising from an
$\alpha$-affine action is therefore an $\alpha$-Lyndon length
function. (Moreover $\mathbf{b}(g)=\ell(g)=\min\{d(x,gx):x\in X\}$
for $g\in\ker\alpha$ if $\mathbf{b}$ arises from a based length
function and $g$ is not an inversion.) We will shortly prove the
converse.

\begin{lemma}\label{length-props} Let $L$ be an $\alpha$-Lyndon length function on $G$.
Then for all $g,h\in G$ we have
\begin{enumerate}
\item $L(g)\geq 0$;
\item $L(gh)\leq L(g)+\alpha_gL(h)$\label{triineq};
\item for $g_1,g_2,\ldots g_n\in G$, if $\bar{g_0}=1$ and $\bar{g_k}=\overline{g_{k-1}}g_k$ for $1\leq k\leq n$ then
$L(g_1g_2\cdots g_n)\leq
\sum_{k=1}^n\alpha_{\overline{g_{k-1}}}L(g_k)$;
\item $\mathbf{c}(g,h)\geq 0$;
\item $\mathbf{c}(g,g)=L(g)$;
\item $\mathbf{c}(g,1)=0$;
\item $\alpha_gL(g^{-1}h)=\alpha_hL(h^{-1}g)$;
\item $\mathbf{c}(g,h)=\mathbf{c}(h,g)\leq\min\{L(g),L(h)\}$;
\end{enumerate}
\end{lemma}
\pf It is clear from the definition of $\mathbf{c}$ and of
$\alpha$-Lyndon length functions that $\mathbf{c}(g,g)=L(g)$ and
$\mathbf{c}(g,1)=0$. One has
$\alpha_hL(h^{-1}g)=\alpha_g\alpha_{g^{-1}h}L(h^{-1}g)=\alpha_gL(g^{-1}h)$
by (L3) whence $\mathbf{c}(g,h)=\mathbf{c}(h,g)$. Taking $k=1$ in
the isosceles condition now gives $\mathbf{c}(g,h)\geq 0$. It
follows that $L(g)=\mathbf{c}(g,g)\geq 0$. Taking $g=k$ in the
isosceles condition gives $\mathbf{c}(g,h)\leq L(g)$; consequently
$\mathbf{c}(g,h)=\mathbf{c}(h,g)\leq L(h)$.

Now $0\leq
2\mathbf{c}(g^{-1},h)=L(g^{-1})+L(h)-\alpha_{g^{-1}}L(gh)$, giving
$L(g)-L(gh)+\alpha_gL(h)=\alpha_g(L(g^{-1})+L(h)-\alpha_{g^{-1}}L(gh))\geq
0$, whence $L(gh)\leq L(g)+\alpha_gL(h)$. The inequality in (3)
follows by induction on $n$.\qed\\

Recall that a $\Lambda$-metric space $(X,d)$ is
\emph{$\delta$-hyperbolic} (with respect to $v$) if $(x\cdot
y)\geq\min\{(x\cdot z),(y\cdot z)\}-\delta$ for $x,y,z\in X$: here
$(x\cdot y)=(x\cdot
y)_{v}=\frac{1}{2}\left(d(x,v)+d(y,v)-d(x,y)\right)$. Indeed, this
`Gromov inner product' gives a convenient characterisation of
$\Lambda$-trees which we will use later: a $\Lambda$-metric space is
a $\Lambda$-tree if and only if, for all $x,y,v\in X$ we have
\be\item $(x\cdot y)_{v}\in\Lambda$; \item $X$ is geodesic; and
\item $X$ is 0-hyperbolic (with respect to  $v$). \ee

In fact it is sufficient to check these conditions for a particular
choice of $v$ while $x$ and $y$ vary.

Three points of a $\Lambda$-tree $X$ are \emph{collinear} if there
is a segment in $X$ containing all three. A $\Lambda$-tree is
\emph{linear} if $x$, $y$ and $z$ are collinear for all $x,y,z\in
X$. Equivalently $X$ is linear if it admits an isometric embedding
in $\Lambda$ (see \cite[\S2.3]{Chiswell-book}).

We write $[x,y,z]$ as a common shorthand for the segment $[x,z]$
where $y\in[x,z]$, and for the assertion that $y\in[x,z]$. More
generally, $[x_1,\ldots,x_n]$ denotes the segment $[x_1,x_n]$ where
$$d(x_1,x_m)=\sum_{k=1}^m d(x_{k-1},x_k)\mbox{\ for all\ }1\leq
m\leq n$$
and the assertion that the displayed equation holds.
If $x_{i-1}\neq x_i$ we may replace the comma between them by a
semicolon. For example $[x,y,z;w]$ implies that $y\in[x,z]$ and
$z\in[x,w]$, and that $z\neq w$. This of course implies that $x\neq
w$ and $y\neq w$ also, and that $z\in[y,w]$.

In general the segments $[x,y]$, $[x,z]$ and $[y,z]$ have exactly
one point in common, which is denoted $\Y(x,y,z)$ (see
\cite[2.1.2]{Chiswell-book}).

\begin{lemma}\label{affine-coll}
\be\item Let $\phi$ be an $\alpha_*$-affine isomorphism $X\to X'$
where $X$ and $X'$ are $\Lambda$-trees. Then
$y\in[x,z]\Leftrightarrow \phi(y)\in[\phi(x),\phi(z)]$ and
$u=\Y(x,y,z)\Leftrightarrow\phi(u)=\Y(\phi(x),\phi(y),\phi(z))$.
\item
Suppose that $\phi:X\to X'$ is a function, $v\in X$ and that for
linear subtrees $J$ of $X$ with $v$ as an endpoint the restriction
$\phi|_J$ is $\alpha_*$-affine. If
$\phi\left(\Y(x,y,v)\right)=\Y(\phi(x),\phi(y),\phi(v))$ then $\phi$
is $\alpha_*$-affine. \ee\end{lemma}

\pf Note that $y\in[x,z]$ if and only if $d(x,y)+d(y,z)=d(x,z)$.
Applying $\alpha_*$ to both sides of this equation gives the
required implications for the first assertion of (1). For the second
assertion, note that $u$ is characterised by the fact that it
belongs to the three segments $[x,y]$, $[x,z]$ and $[y,z]$. By the
first assertion, $\phi(u)$ belongs to the images of these segments
under $\phi$, giving the result.

To establish (2), for $x,y\in X$, we put $u=\Y(x,y,v)$. Then
\begin{eqnarray*}d(\phi(x),\phi(y))&=&d(\phi(x),\phi(u))+d(\phi(u),\phi(y))\\
&=&\alpha_*d(x,u)+\alpha_*d(u,y)\\ &=&\alpha_*d(x,y).
\end{eqnarray*}
\qed\\

\begin{theorem}\label{basechange}
\begin{enumerate}\item Let $(X,d)$ be a $\Lambda$-metric space and
$v\in X$ such that
\begin{enumerate}\item $(x\cdot y)_{v}\in \Lambda$; \item
$X$ is 0-hyperbolic with respect to $v$.\end{enumerate} Then there
exists a $\Lambda$-tree $(X',d')$ and an isometric embedding
$\phi:X\to X'$ such that if $\psi:X\to Z$ is any $\alpha_*$-affine
map (where $Z$ is a $\Lambda$-tree) then there is a unique
$\alpha_*$-affine map $\mu:X'\to Z$ such that $\mu\circ\phi=\psi$.
\item In the notation of part (1), an $\alpha$-affine action on the $\Lambda$-metric space $X$ has a unique extension to an $\alpha$-affine action on the $\Lambda$-tree $X'$.

\item Let $(X_1,d_1)$ be a $\Lambda_1$-tree,
$\alpha:G\to\Aut^+(\Lambda_1)$ a homomorphism, let an
$\alpha$-affine action of $G$ on $X_1$ be given and let
$h:\Lambda_1\to\Lambda_2$ be an order preserving homomorphism.
Suppose that $\bar{\alpha}:G\to\Aut^+(\Lambda_2)$ is a homomorphism
such that $h\circ\alpha_g=\bar{\alpha_g}\circ h$ for all $g\in G$.
Then there is a $\Lambda_2$-tree $(X_2,d_2)$, an
$\bar{\alpha}$-affine action of $G$ on $X_2$ and a map $\phi:X_1\to
X_2$ satisfying\begin{enumerate}
\item $d_2(\phi(x),\phi(y))=h\
d_1(x,y)$ for $x,y\in X_1$;
\item $\phi(gx)=g\phi(x)$ for $x\in
X_1$.
\item for $p\in X_1$, we have $L_2(g)=h\ L_1(g)$ where $L_1$ denotes the based length function with respect to $p$ and $L_2$ the based length function with respect to $\phi(p)$.
\end{enumerate}
\end{enumerate}
\end{theorem}

\pf (1) The $\Lambda'$-tree $(X',d')$ is as constructed in
\cite[Theorem 2.4.4]{Chiswell-book}. We will follow the notation
used there. To prove our assertion, it suffices to show that a
suitable map $\mu$ can be found for a given $\psi$. Obviously such a
map must satisfy $\mu\langle x,d(x,v)\rangle=\mu(\phi(x))=\psi(x)$.
Moreover, for $x\in X$ and $m\geq 0$ in $\Lambda$ with $m\leq
d(x,v)$, if we let $x_m$ be the point on the segment
$[\psi(v),\psi(x)]_{X'}$ at distance $\alpha_*m$ from $x_0=\psi(v)$,
then any $\mu$ with the required properties must satisfy $\mu\langle
x,m\rangle=x_m$. This amounts to a definition of $\mu$, and it
remains only to show that $\mu$ is $\alpha_*$-affine. Now for $n\leq
m\leq d(x,v)$ we have
\begin{eqnarray*}d_Z(\mu\langle x,m\rangle,\mu\langle
x,n\rangle)&=&d_Z(x_m,x_n)\\
&=& d_Z(x_m,x_0)-d_Z(x_n,x_0)\\
&=&\alpha_*m-\alpha_*n\\
&=&\alpha_*(m-n)\\ &=& \alpha_*d_{X'}(\langle x,m\rangle,\langle
x,n\rangle).\end{eqnarray*} So  $\mu$ restricted to subtrees $J$
with $v$ as an endpoint is therefore $\alpha_*$-affine. The result
now follows from Lemma~\ref{affine-coll}(2).

(2) Fix an $\alpha$-affine action of $G$ on $X$. Let $Z=X'$, for
$g\in G$ let $\psi=\phi\circ g$, and let $\mu_g$ be the associated
map $X'\to X'$. We claim that the assignation $g\cdot
x_m=\mu_g(x_m)$ defines the required $\alpha$-affine action of $G$
on $X'$. That each $\mu_g$ is $\alpha_g$-affine follows from part
(1).

Now note that $\mu_{g}\circ\phi=\phi\circ g$, whence
\begin{eqnarray*}\left(\mu_{g}\circ \mu_h\right)\circ\phi &=& \mu_{g}\circ\phi\circ h\\
&=& \phi\circ\left(g\circ h\right),\end{eqnarray*}

So $\mu_{gh}$ and $\mu_g\circ\mu_h$ both satisfy
$\mu\circ\phi=\phi\circ\left( g\circ h\right)$. Since such a map is
unique by part (1), we find that $\mu_{gh}=\mu_g\circ\mu_h$. It is
also easy to see that $\mu_1=\id_{X'}$, giving
$\id_{X'}=\mu_{g^{-1}}\circ\mu_g=\mu_g\circ\mu_{g^{-1}}$. This shows
that $\mu_g$ and $\mu_{g^{-1}}$ are mutually inverse functions on
$X'$. In particular, they are bijections.

(3) The $\Lambda_2$-tree $(X_2,d_2)$ is the the base change functor
as described in \cite[2.4.7]{Chiswell-book}, where property (a) is
shown.
A 0-hyperbolic $\Lambda_2$-metric space $Z$ is constructed, whose
points are equivalence classes $\langle x\rangle$ of points of $X_1$
under the relation $x\sim y$ if $h\ d_1(x,y)=0$. We put
$\phi(x)=\langle x\rangle$. It remains to show that the natural
action of $G$ on $X_2$ ($g\langle x\rangle=\langle gx\rangle$) is
$\bar{\alpha}$-affine. (Part (b) is then clear, and the last part is
routine.) So observe that
\begin{eqnarray*}d_2(\langle gx\rangle,\langle gy\rangle)&=h\
d_1(gx,gy)\\ &=h\ \alpha_g d_1(x,y)\\ &=\bar{\alpha}_g h\ d_1(x,y)\\
&=\bar{\alpha}_g d_2(\langle x\rangle,\langle
y\rangle).\end{eqnarray*}\nopagebreak\qed

The $\Lambda_2$-tree $X_2$ of Theorem~\ref{basechange}(3) is called
the \emph{base change functor} and is denoted
$\Lambda_2\otimes_{\Lambda_1}X_1$. (Note however that it depends on
$h$.)

The existence of a homomorphism $\bar{\alpha}$ as in
Theorem~\ref{basechange} implies that $\ker(h)$ is
$\alpha_g$-invariant for all $g$. If $h$ is surjective, this
condition is necessary and sufficient. Otherwise, this condition
ensures that a suitable $h(\Lambda_1)$-tree can be found, and the
problem reduces to extending the automorphisms $\bar{\alpha_g}$ of
$h(\Lambda_1)$ to automorphisms of $\Lambda_2$. However extending
automorphisms in this way is not possible in general: for example,
$\Aut^+(\langle\mathbb{Q},\sqrt{2}\rangle)$ is trivial, while
$\Aut^+(\mathbb{Q})$ is not.

As in the isometric case, an important special case arises when we
take $\Lambda_1=\Lambda_2$ and $h(\lambda)=2\lambda$. The resulting
$\Lambda$-tree is called the \emph{barycentric subdivision} of $X$.
Taking $\bar{\alpha}=\alpha$, it is easy to see that the hypotheses
of the theorem are satisfied, so an affine action on $X$ extends to
an affine action on the barycentric subdivision in general. As in
the isometric case, if a segment $[x,y]$ is stabilised by $g$ then
there is a point of the barycentric subdivision fixed by $g$, namely
the midpoint of the segment.

\begin{theorem}\label{based-lyndon} Let $\Lambda$ be an ordered abelian group, $G$ a group and
$\alpha:G\to\Aut^+(\Lambda)$ a homomorphism.

If $L$ is a $\Lambda$-valued $\alpha$-Lyndon length function on $G$
there exists a $\Lambda$-tree $X$ on which $G$ has an
$\alpha$-affine action that induces $L$.
\end{theorem}

\pf Suppose that a $\Lambda$-valued $\alpha$-Lyndon length function
$L$ is given.
Define $\delta:G\times G\to\Lambda$ by the rule
$\delta(g,h)=\alpha_gL(g^{-1}h)$. It follows from
Lemma~\ref{length-props}(7) that $\delta$ is symmetric. Moreover, by
Lemma~\ref{length-props}(\ref{triineq}) we have $L(g^{-1}k)\leq
L(g^{-1}h)+\alpha_{g^{-1}h}L(h^{-1}k)$. Applying $\alpha_g$ to both
sides of this inequality yields
$\delta(g,k)\leq\delta(g,h)+\delta(h,k)$. We also have
$\delta(\gamma g,\gamma h)=\alpha_\gamma\delta(g,h)$ and
$\delta(g,g)=0$, whence $(G,\delta)$ is a pseudometric space.
Denoting the corresponding metric space by $(\bar{G},\bar{\delta})$,
we also see that right multiplication by $G$ gives an
$\alpha$-affine action. Taking $\bar{1}$ to be the basepoint for the
Gromov inner product on $\bar{G}$, it is routine to show that
$(\bar{g}\cdot\bar{h})_{\bar{1}}=\c({g},{h})$. Thus
$(\bar{g}\cdot\bar{h})_{\bar{1}}\in\Lambda$.

To show that $(\bar{G},\bar{\delta})$ is a 0-hyperbolic
$\Lambda_0$-metric space with respect to $\bar{1}$, it suffices to
show that for all $\bar{g},\bar{h},\bar{k}\in\bar{G}$ we have the
implication
\begin{eqnarray*} & \bar{\delta}(\bar{g},\bar{1}) +
\bar{\delta}(\bar{h},\bar{1})-\bar{\delta}(\bar{g},\bar{h})
<\bar{\delta}(\bar{g},\bar{1})+\bar{\delta}(\bar{k},\bar{1})-\bar{\delta}(\bar{g},\bar{k})\\
\Rightarrow &\bar{\delta}(\bar{g},\bar{1}) +
\bar{\delta}(\bar{h},\bar{1})-\bar{\delta}(\bar{g},\bar{h})=
\bar{\delta}(\bar{h},\bar{1})+\bar{\delta}(\bar{k},\bar{1})
-\bar{\delta}(\bar{h},\bar{k}). \end{eqnarray*} This follows from
the isosceles condition on $\mathbf{c}$.

So $G$ has an $\alpha$-affine action on the 0-hyperbolic
$\Lambda_0$-metric space $\bar{G}$ with
$\bar{\delta}(\gamma\bar{1},\bar{1})=L(\gamma)$. By
Theorem~\ref{basechange}, $\bar{G}$ has a $G$-equivariant embedding
in a $\Lambda$-tree with the required properties.\qed

As in the isometric case, the $\Lambda$-tree $X$ in
Theorem~\ref{based-lyndon} is unique in the following sense. If
$(Z,d'')$ is a $\Lambda$-tree on which $G$ has an $\alpha$-affine
action, and $w\in Z$ satisfies $L=L_w$ then there is a unique
$G$-equivariant isometry $\mu:X\to Z$ such that $\mu(x)=w$ and the
image of $\mu$ coincides with the subtree of $Z$ spanned by the
orbit $Gw$. The proof of this assertion closely follows the proof in
the isometric case (see \cite[Theorem 4.6]{Chiswell-book}) and will
be omitted.

In the literature a group is said to be $\Lambda$-free if it admits
a free isometric action (without inversions) on a $\Lambda$-tree,
and tree-free if it is $\Lambda$-free for some $\Lambda$. We will
say that a group is ITF (isometric tree-free) if it admits a free
isometric action without inversions on a $\Lambda$-tree for some
$\Lambda$. A group is ATF (affine tree-free) if it admits a free
affine action without inversions on a $\Lambda$-tree for some
$\Lambda$. If we want to specify the ordered abelian group\ we will
refer instead to an $\ITF(\Lambda)$ or $\ATF(\Lambda)$ group as
appropriate. Note that henceforth a free action will be assumed to
be without inversions unless the contrary is explicitly stated.

\subsection{Affine automorphisms of $\Lambda$-trees}

Let $g$ be an affine automorphism of a $\Lambda$-tree $X$. Let $X^g$
denote the set of points of $X$ fixed by $g$. Note that if
$\alpha_g\neq 1=\id_\Lambda$ the set of all fixed points need
\emph{not} form a subtree of $X$, as the map
$(x,y,z)\mapsto(x,y,z+y+x)$ in the case $X=\Lambda=\Z^3$ shows. Let
$A_g=\{x\in X: x\in[g^{-1}x,gx]\}$, and $\tilde{A_g}=\{x\in X:
g^{-1}x,x,gx\mbox{ are collinear}\}$. Call $g$ \emph{elliptic} if
$X^{g}\neq\emptyset$. Note that $X^g\subseteq X^{g^n}$ and
$A_g\subseteq A_{g^n}$ for $n\neq 0$, and $A_g\subseteq
\tilde{A}_g$.

We call $g$ an \emph{inversion} if $X^g=\emptyset$ but
$X^{g^2}\neq\emptyset$. In this case there exists a segment $[x,y]$
in $X$ whose endpoints are swapped by $g$, and $g^2$ stabilises the
segment. Note that like in the isometric case, if $g $ is an
inversion of $X$ then $g$ has a unique fixed point in the
barycentric subdivision.

Call $g$ a \emph{nesting reflection} if $X^{g^2}=\emptyset$ and
there exists $x\in X$ such that $[x;gx;g^{-1} x]$ or $[x;g^{-1}x;g
x]$. In this case $g^n$ fixes no point for $n\neq 0$, and either
$g[x,g x]\subset[x,gx]$ or $g[x,g x]\supset[x,gx]$. See
Example~\ref{nesting-ex} for an example of this behaviour.

If $g$ does not satisfy any of the criteria above, then $g$ is said
to be \emph{hyperbolic}. As we will see (Theorem~\ref{hyper}), a
hyperbolic affine automorphism behaves much like its isometric
counterpart.

Let us clarify the situation that arises if there is some point
$x\in\tilde{A}_g\backslash A_g$ --- a situation that cannot arise in
the isometric case unless $g$ is an inversion.

Suppose that $g$ is an affine automorphism of $X$, and that $x\in X$
satisfies $[x;gx,g^{-1}x]$ or $[x;g^{-1}x,gx]$. We lose no
generality in the following discussion by assuming the former.
Suppose further that $y\in A_g$ --- thus $y\in X$ satisfies
$[g^{-1}y,y,gy]$. Suppose initially that $y$ satisfies
$[y,g^nx,g^{n+1}x]$ for some $n\in\Z$. Then replacing $y$ by $g^{\pm
1}y$ if necessary, we can ensure that
$[g^{-1}y,y,gy,g^nx;g^{n+1}x,g^{n-1}x]$ or $[gy,y,g^{-1}y,
g^nx;g^{n+1}x,g^{n-1}x]$. But in either case we get a contradiction:
in the first case, we have $g[g^nx,y]\not\subseteq[g^nx,y]$
contradicting $g^{-1}[g^nx,y]\supset [g^nx,y]$; while in the second
case, we obtain $g[g^nx,y]\supset[g^nx,y]$ which contradicts
$g^{-1}[g^nx,y]\not\subseteq[g^nx,y]$.

Therefore for any $y\in X$ satisfying $[g^{-1}y,y,gy]$ we must have
$v=\Y(y,g^nx,g^{n+1}x)\notin\{g^nx,g^{n+1}x\}$ for all $n\in\Z$ and
thus $v\in I=\cap_{n\in\Z}[g^nx,g^{n+1}x]$. If $y\in I$ (so $y=v$)
and $gy\neq y$ then (replacing $x$ by $gx$ if necessary) we have
$[g^nx,g^{n+2}x,\ldots, g^{-1}y;y;gy,\ldots
g^{n+1}x,g^{n-1}x]\subseteq I$, which forces a contradiction; thus
$gy=y$. This shows that $A_g\cap I$ consists of points fixed by $g$
(though it may well be empty).

Let us show next that $g$ can fix at most one point of $[x,gx]$. For
if $p,q\in [x,gx]$ are both fixed by $g$ then $p,q\in I$ and we have
$[x,g^2x,\ldots, p,q,\ldots, g^3x,gx]$ which implies
$[g^{-1}x,gx,\ldots,g^{-1}p,g^{-1}q,\ldots,g^2x,x]$. Thus we have
both $[x,p,q,gx]$ and $[x,q,p,gx]$, giving $q=p$.

Finally, note that for any $y$ such that $[g^{-1}y,y,gy]$ holds, if
$w$ belongs to the segment $[g^{-1}y,y,gy]$ and to $[x,gx]$ then
$[g^{-1}w,w,gw]$ holds also, so we can replace $y$ by $w$ in the
argument above to conclude that $gw=w$, and hence $g^{-1}y=y=gy$.
Otherwise $[g^{-1}y,y,gy]$ has trivial intersection with $[x,gx]$,
and $v=\Y(g^\epsilon y,g^nx,g^{n+1}x)$ for $\epsilon\in\{-1,0,1\}$.
We claim that $v$ is fixed by $g$ (and is thus the unique fixed
point of $g$ in the segment $[x,gx]$). For
$[g^nx;g^{n+1}x,g^{n-1}x]$ holds, which implies
\begin{eqnarray*}v&=&\Y(g^{-1}y,g^{n-1}x,g^nx)\\ &=&
g^{-1}\Y(y,g^nx,g^{n+1}x)\\ &=& g^{-1}v.\end{eqnarray*}\\

In summary,
\begin{lemma}\label{nesting-lemma}
Suppose that $[x;gx,g^{-1}x]$ or $[x;g^{-1}x,gx]$ holds. \be\item
Then $g$ has at most one fixed point $v$ in $[x,gx]$, in which case
$v\in\cap_{n\in\Z}[g^nx,g^{n+1}x]$ and any other fixed point $y$ in
$X$ must satisfy $v=\Y(y,x,gx)$.
\item
If $y\in A_g$ then $\Y(y,gx,x)$ is a fixed point of $g$. If $g$ is
not elliptic then $A_g=\emptyset$.
\item If $g^2$ is not elliptic
then $g$ is not an inversion; thus $g$ is a nesting reflection.
\item If $h$ is an affine automorphism of $X$ with $A_h\neq\emptyset$ and $h$ fixes no point
then
$h$ is hyperbolic.
\item If $h$ is a nesting reflection then
$h^2$ is hyperbolic.
\item If $h$ is hyperbolic, then
$h^n$ is hyperbolic for all $n\neq 0$.\qed\ee
\end{lemma}

The case where $g$ is hyperbolic is described further in the next
theorem. For now let us warn the reader that an affine automorphism
may behave near one point like an elliptic automorphism, and
elsewhere like a hyperbolic automorphism. Moreover, unlike the
isometric case, it is not always possible to distinguish between the
cases we describe by means of based length functions.

\example\label{mixed-ex} Fix a real number $0<a<1$ and let $X$ be
the quotient of $[0,\infty)_{\R}\times\{1,2,3\}$, where $(0,i)$ is
identified with $(0,j)$ for all $i$ and $j$, and put
$$d((x,i),(y,j))=\left\{\begin{array}{cl}|x-y| & \mbox{if }i=j\\
|x|+|y| & \mbox{otherwise}\end{array}\right.$$ Define $g:X\to X$ via
$$g:(x,i)\mapsto \left\{\begin{array}{cl}(ax,\ 3-i) &
i\in\{1,2\}\\ (ax,\ 3) & i=3.
\end{array}\right.$$

Then $g$ is $\alpha_*$-affine where $\alpha_*x=ax$. Since (the
equivalence class of) $(0,1)$ is fixed by $g$, the latter is clearly
elliptic. However, taking $X_3$ to be the subtree consisting of
points of the form $(x,3)$ with $x\neq 0$, we see that $X_3$ is
$g$-invariant and the restriction of $g$ to this subtree is
hyperbolic.

On the other hand, taking $X_{12}$ to be the complement of $X_3$, we
see that this is also a $g$-invariant subtree and that
$[p;gp;g^{-1}p]$ holds for all $p\neq (0,i)$. Thus (denoting the
restriction of $g$ to $X_{12}$ also by $g$) we have
$\tilde{A}_g=X_{12}$ while $A_g$ consists of a single point, which
is fixed by $g$. It is easy to show that $\b_p(g)<0$ for all $p\in
X_{12}$ where $p$ is not the fixed point, while $\b_p(g)>0$ for
$p\in X_3$. Of course $\b_p(g)=0$ precisely when $p$ is the fixed
point. (See \S1.2 for the definition of $\b$.)

\example\label{nesting-ex} Let us also describe an example of a
nesting reflection. Let $X=\Lambda=\R$. Put
$g\cdot x=\displaystyle{2-\frac{x}{2}}$. Then the unique fixed point
of $g^2$ is $4/3$, which is also the unique fixed point of $g$.
However if we replace $\Lambda=\R$ by the subgroup $\Lambda_0$
consisting of the dyadic rationals ($a/2^n$), then $g$ stabilises
$\Lambda_0$, but $g$ and $g^2$ have no fixed point,
since $4/3$ is not an element of $\Lambda_0$.\\

However, if automorphisms are assumed to be rigid, the situation
becomes easier to analyse. An automorphism $g$ of a $\Lambda$-tree
is \emph{rigid} (or \emph{non-nesting}) if no closed segment is
mapped properly into itself by $g$ or $g^{-1}$. In this case,
nesting reflections described above are not possible (as the
nomenclature suggests), and if $g$ fixes a point, the set of all
fixed points forms a subtree. If $g$ is rigid and $g^2$ does not fix
a point, then $g$ is hyperbolic in a similar sense to the isometric
case: there is a maximal linear $g$-invariant subtree $A_g$, such
that $x<g x$ for all $x\in A_g$ (with respect to one of the natural
linear orders on $A_g$).

Note that if $g$ is hyperbolic and $A_g$ is spanned by a single
$\langle g\rangle$-orbit then $g$ is rigid. (Such a $g$ is said to
be archimedean.)

In \cite{Hudson}, Hudson considers rigid group actions on pretrees,
and shows that much of the theory of isometries and of isometric
actions can be extended to this general setting.

Since the rigid case is better behaved than the general case, we
will focus much of our attention on the rigid case, and show where
possible that rigidity is preserved by our constructions. We make no
such restriction in the next theorem however. First, a lemma.

\begin{lemma}\label{coll}
If $I_1$ and $I_2$ are disjoint linear subtrees of $X$, then either
$I_1\cup I_2$ is collinear or, for $i=1$ or $i=2$, there exists
$x_0\in I_i$ with $[x_0,z]\cap I_{i}=\{x_0\}$ for all $z\in
I_{3-i}$. The point $x_0$ belongs to every segment joining a point
of $I_1$ to a point of $I_2$.
\end{lemma}

\pf Assume that $I_1\cup I_2$ is not collinear. Then there exist
$x,y,z\in I_1\cup I_2$ with $u=\Y(x,y,z)\notin\{x,y,z\}$. Without
loss of generality, $x,y\in I_1$ and $z\in I_2$.
We claim that $i=1$ and $x_0=u$ give the required properties.
Firstly, suppose that $v\in[u,z]\cap I_1$. Then we have $[z;v,u,x]$
and $[z;v,u,y]$, so that $u=\Y(x,y,z)=\Y(x,y,v)$. If $u\neq v$ this
contradicts the linearity of $I_1$, since $x,y,v\in I_1$. This shows
that $[u,z]\cap I_1=\{u\}$. Now if $z'$ is any point of $I_2$, we
have \begin{eqnarray*}[u,z']\cap I_1 &\subseteq & ([u,z]\cap I_1)
\cup ([z,z']\cap I_1)\\ &=& \{u\} \cup \emptyset.\end{eqnarray*}

Now suppose that $s\in I_1$ and $t\in I_2$; we claim that
$x_0=u\in[s,t]$. Put $w=\Y(s,t,u)$ and observe that $w\in I_1$,
since $s,u\in I_1$. Thus
\begin{eqnarray*}w&\in&[w,t]\cap
I_1\\ &\subseteq &[u,t]\cap I_1\\
&=&\{u\},\end{eqnarray*} forcing $u=w\in[s,t]$. \qed\\

Note that for any affine automorphism $g$ we have the implications
$$\begin{array}{rcl}u\in A_g&\Leftrightarrow& [g^{-1}u,u,gu]\\
&\Leftrightarrow& [\gamma g^{-1}u,\gamma u, \gamma gu]\\
&\Leftrightarrow& [(\gamma g^{-1}\gamma^{-1}) (\gamma u),\gamma u,
(\gamma g\gamma^{-1}) (\gamma u)\\ &\Leftrightarrow& \gamma u\in
A_{\gamma g\gamma^{-1}}
\end{array}$$
Thus $\gamma\cdot A_g=A_{\gamma g\gamma^{-1}}$. Similarly
$\gamma\cdot \tilde{A_g}=\tilde{A}_{\gamma g\gamma^{-1}}$.

If $I$ and $J$ are disjoint subtrees of $X$, we use the notation
$\overline{\Br}(I,J)$ for the (closed) bridge between $I$ and $J$
(see \cite[\S2.1]{Chiswell-book}). We will also write $\Br^o(I,J)$
for the open bridge between $I$ and $J$: this is the set
$[x,y]\backslash(I\cup J)$ where $x\in I$ and $y\in J$, and can be
shown to be a linear subtree independent of these points.

The following theorem follows closely the argument given in
\cite[3.1.4]{Chiswell-book} for the isometric case, but is in fact
valid in the general situation where $X$ is a median pretree and $g$
is any pretree automorphism.
(except of course for the last two parts which refer to $d$).

\begin{theorem}\label{hyper}
Suppose that $g$ is hyperbolic (that is, suppose that $g$ fixes no
point of $X$ and is not a nesting reflection or an inversion). For
$x\in X$, we put $u=u_x=u(x,g)=\Y(g^{-1}x,x,gx)$.

Then
\begin{enumerate}
\item $A_g$ is non-empty;
\item $A_g$ is linear;
\item $A_g$ is $\langle g\rangle$-invariant;
\item $A_g$ is a closed subtree;
\item $A_g$ is not properly contained in a linear subtree of $X$;
\item If $T$ is a maximal linear $g$-invariant subtree of $X$, then
$L=A_g$;
\item If $x\in X$ then $[x,u]=\overline{\Br}(x,A_g)$;
\item $[x,gx]\cap A_g=[u,gu]$;
\item $[x,gx]=[x,u,gu,gx]$;
\item $A_g=A_{g^n}$ for $n\neq 0$;
\item $\a_x(g)=d(x,u)$;
\item $\b_x(g)=d(u,gu)$.

\end{enumerate}
\end{theorem}

\pf Let $x\in X$, and let $u=u_x$ be as in the Theorem.
Now $u\in[g^{-1}x,x]$ implies $gu\in[x,gx]$. Also $u\in[x,gx]$, so
we have either $[x,u,gu,gx]$ or $[x,gu,u,gx]$. But in the latter
case we have $u\in[gu,gx]$ and $gu\in[u,x]$, forcing
$u,g^2u\in[gu,gx]$. Thus either $[gu,u,g^2u]$ holds, giving
$[u,g^{-1}u,gu]$ or $[gu,g^2u,u]$ holds, giving $[u,gu,g^{-1}u]$.
Both of these cases force $g$ to be either a nesting reflection, an
inversion or elliptic, contradicting our original assumption.

Therefore $[x,u,gu,gx]$  holds, giving part (9). We also have
$[g^{-1}x,g^{-1}u, u,x]$. Now $u\in[g^{-1}x,gx]$ implies
$[g^{-1}x,g^{-1}u, u,gx]$, whence $[g^{-1}x,g^{-1}u, u,gu,gx]$ and
$[g^{-1}u,u, gu]$. Therefore $u\in A_g$. This shows that
$A_g\neq\emptyset$, establishing part (1).

Next, for $p\in A_g$ put $A_p=\cup_{n\in\Z}[g^np,g^{n+1}p]$. Then,
since $[g^np;g^{n+1}p;g^{n+2}p]$ holds for $n\in\Z$, the Piecewise
Geodesic Proposition (see \cite[2.1.5]{Chiswell-book}) and an
induction argument give $[g^np,g^mp,g^lp]$ for $n\leq m\leq l$,
whence $A_p$ is linear.

Now if $q\in A_p$, we have $[g^np,q,g^{n+1}p]$ for some $n$, whence
$[g^{n-1}p,g^{-1}q,g^{n}p]$ and $[g^{n+1}p,gq,g^{n+2}p]$, which
gives $[g^{n-1}p,g^{-1}q,g^{n}p,q,g^{n+1}p,gq,g^{n+2}p]$, and hence
$[g^{-1}q,q,gq]$. This gives $q\in A_g$. Therefore $A_p\subseteq
A_g$, and consequently $\cup_{p\in A_g}A_p\subseteq A_g$. It is now
easy to see that $A_g$ is $\langle g\rangle$-invariant, since each
$A_p$ clearly is. This proves (3)

For $p\in A_g$, put $A_p^+=\cup_{n\geq 0}[p,g^np]$, and
$A_p^-=\cup_{n\leq 0}[p,g^np]$. Then $A_p=A_p^+\cup A_p^-$ and
$A_p^+\cap A_p^-=\{p\}$. We observe that if $J$ is a subtree of
$A_p$ containing $p$ then  $gJ\subseteq J$ implies $A_p^+\subseteq
J$, and $J\subseteq gJ$ implies $A_p^-\subseteq J$.

\begin{claim}
$A_p$ is a minimal (non-empty) $\langle g\rangle$-invariant subtree
of $X$. Moreover, for $p,q\in A_g$,  either $A_p=A_q$ or $A_p\cap
A_q=\emptyset$.
\end{claim}
To see this, suppose that $J$ is a non-empty $\langle
g\rangle$-invariant subtree of $A_p$. If $x\in J$ then
$x\in[g^mp,g^{m+1}p]$ for some $m$. Thus
$g^{-1}x\in[g^{m-1}p,g^mp]$, giving
$[g^{m-1}p,g^{-1}x,g^mp,x,g^{m+1}p]$, whence
$g^mp\in[g^{-1}x,x]\subseteq J$. Since $J$ is $\langle
g\rangle$-invariant, this gives $p\in J$. Applying the observation
above to $J$, we obtain $A_p^+\subseteq J$ and $A_p^-\subseteq J$,
that is, $A_p\subseteq J$. This forces $J=A_p$.

For the second assertion, put $J=A_p\cap A_q$. If $J$ is non-empty,
then it is a $\langle g\rangle$-invariant subtree of $A_p$, which
now implies $J=A_p$, whence $A_p=A_q$.

\begin{claim}
$A_p\cup A_q$ is linear for $p,q\in A_g$.
\end{claim}

Suppose that $A_p\cup A_q$ is not linear. Then $A_p$ and $A_q$ are
disjoint and linear, so swapping $p$ and $q$ if necessary, we see
that there exists $p_0\in A_p$ with $[p_0,q']\cap A_p=\{p_0\}$ for
all $q'\in A_q$, by Lemma~\ref{coll}. Moreover $p_0\in[p',q']$ for
all $p'\in A_p$ and $q'\in A_q$.

But $[gp_0,gq']$ is a segment joining a point of $A_p$ to a point of
$A_q$, forcing $p_0\in[gp_0,gq']\cap A_p=\{gp_0\}$. Thus $gp_0=p_0$,
a contradiction.

\begin{claim}
Any three points of $A_g$ are collinear.
\end{claim}

Choose $p,q,r\in A_g$. If $A_p$, $A_q$ and $A_r$ are not pairwise
disjoint, then two of these subtrees coincide, and this claim
reduces to the previous claim. So suppose that they are pairwise
disjoint. If $p$, $q$ and $r$ are not collinear, let $w=\Y(p,q,r)$.
Then $w\in[p,q]$. If $w\in A_p$, then $w\in[q,r]$, giving
$g^nw\in[g^nq,g^nr]$, a segment with endpoints belonging to $A_q$
and $A_r$. But the subtree $I$ spanned by $A_q\cup A_r$ is linear,
with $A_q$ and $A_r$ as $\langle g\rangle$-invariant subtrees, so
this forces $A_p=A_w=\cup_{n\in\Z}[g^nw,g^{n+1}w]\subseteq I$,
whence $A_p$, $A_q$ and $A_r$ are collinear, a contradiction.

So we may suppose that $w\notin A_p$, and likewise that $w\notin
A_q\cup A_r$. Then $w\in\Br^o(A_p,A_q)$ and similarly
$w\in\Br^o(A_p,A_r)\cap\Br^o(A_q,A_r)$, so that $w$ is the unique
point common to all three open bridges. By $g$-invariance of $A_p$,
$A_q$ and $A_r$, we also see that $gw$ is the unique point common to
the three open bridges, forcing $gw=w$, a contradiction. This proves
part (2).

\begin{claim}
$A_g$ is a subtree.
\end{claim}
Let $p,q\in A_g$; it suffices to show that $[p,q]\subseteq A_g$. If
$A_p\cap A_q\neq\emptyset$, then $A_p=A_q$. It is clear in this case
that $[p,q]\subseteq A_p\subseteq A_g$, since $A_p$ is a subtree.

So suppose that $A_p\cap A_q=\emptyset$, and let
$r\in\Br^o(A_p,A_q)$. Then $gr\in\Br^o(gA_p,gA_q)=\Br^o(A_p,A_q)$,
and likewise $g^{-1}r\in\Br^o(A_p,A_q)$. Since the bridge between
disjoint subtrees is linear, the points $g^{-1}r$, $r$ and $gr$ are
collinear. Since $g$ is assumed not to be a nesting reflection, this
forces $[g^{-1}r,r,gr]$. That is, $r\in A_g$. This shows that
$\Br^o(A_p,A_q)\subseteq A_g$. Now $[p,q]\subseteq
A_p\cup\Br^o(A_p,A_q)\cup A_q\subseteq A_g$.
\begin{claim}
$A_g$ is a maximal linear subtree of $X$.
\end{claim}

Suppose that $A_g$ is properly contained in a linear subtree $L$ of
$X$ with $x\in L\backslash A_g$. We have shown above that
$u=\Y(g^{-1},x,gx)\in A_g$, and consequently $gu,g^{-1}u\in A_g$.
Swapping the roles of $g$ and $g^{-1}$ if necessary, we now have
$[g^{-1}u,u,gu,x]$, since $A_g$ is a subtree and $x\notin A_g$. But
in general $[x,gx]=[x,u;gu,gx]$, contradicting the configuration
just stated. Therefore $A_g$ is not properly contained in a linear
subtree. This shows that (5) is satisfied. Moreover, since
$A_g\subseteq A_{g^n}$ and $g^n$ is hyperbolic for $n\neq 0$ by
Lemma~\ref{nesting-lemma}, part (10) follows.

\begin{claim}
$A_g$ is a closed subtree of $X$.
\end{claim}

We have seen that $[x,gx]=[x,u_x,gu_x,gx]$, that $u_x,gu_x\in A_g$
and that $A_g$ is a subtree. If $y\in[x,u_x]$ then
$[x,y,u_x,gu_x,gy,gx]$, whence $u_x\in[y,gy]$. Also, $gu_x\in[y,gy]$
giving $u_x\in[g^{-1}y,y]$, and $[g^{-1}x,g^{-1}u_x,u_x,gu_x,gx]$,
giving $u_x\in[g^{-1}y,gy]$. So $u_x=\Y(g^{-1}y,y,gy)=u_y$, giving
$[g^{-1}y,y]\cap[y,gy]=[y,u_x]$. Therefore $y\in A_g$ if and only if
$y=u_x$, giving $[x,u_x]\cap A_g=\{u_x\}$, whence $[gx,gu_x]\cap
A_g=\{gu_x\}$  and $[x,gx]\cap A_g=[u_x,gu_x]$. This proves (8).

For arbitrary $x,y\in X$ with $[x,y]\cap A_g\neq\emptyset$, let
$q\in [x,y]\cap A_g$. Then $[x,q]=[x,u_x,q]$ and $[q,y]=[q,u_y,y]$,
giving $[x,y]=[x,u_x,q,u_y,y]$, and thus $[x,y]\cap A_g=[u_x,u_y]$.
Therefore $A_g$ is a closed subtree, as claimed in (4).\medskip

Writing $u_x=u(x,g)$, we have also shown that
$[x,u(x,g)]=\overline{\Br}(x,A_g)$ where $u(x,g)=\Y(g^{-1}x,x,gx)$,
which establishes (7). Thus $d(x,A_g)=d(x,u(x,g))$.

Now suppose that $T$ is a maximal linear subtree which is
$g$-invariant, and that $x\in T\backslash A_g$.
Since $T$ is $g$-invariant, we have $g^{-1}x, x,gx\in T$. Since
$A_g$ is also $g$-invariant, these points do not belong to $A_g$.
But $u=\Y(g^{-1}x, x,gx)\in\{g^{-1}x, x,gx\}$ since these points
belong to $T$ which is linear, and $u\in A_g$, a contradiction.
Therefore $T\subseteq A_g$, whence (6).

All the assertions have now been established with the exception of
the last two. Note that since $A_g=A_{g^n}$ for $n\neq 0$, we obtain
$u(x,g^n)=u(x,g)$ for $n\neq 0$. It follows that
$L(g^2)=d(x,u)+d(u,gu)+d(gu,g^2u)+d(g^2u,g^2x)$, and
$L(g^{-1})=d(x,u)+d(u,g^{-1}u)+d(g^{-1}u,g^{-1}x)$. Expanding the
expression for $\b_x(g)$, it is now straightforward to show that it
is equal to $d(u,gu)$, and similarly that
$\a_x(g)=d(x,u)$. \qed\\

This theorem shows that in many respects affine hyperbolic
automorphisms behave like their isometric counterparts.

We next give a partial characterisation of the different types of
affine automorphism of a $\Lambda$-tree in terms of based length
functions, via the ancillary function $\b=\b_x$. In fact $\b$ plays
a role similar to the hyperbolic length function $\ell$ in the
isometric case. If $g$ is an isometry and not an inversion then
$\b(g)=\ell(g)$. In particular $\b$ is independent of the basepoint
in this case; there is however no way of escaping the dependence of
$\b$ on the choice of basepoint in general.

\begin{proposition}\label{class-autos}
Let $g$ be an affine automorphism of a $\Lambda$-tree $X$. Then
exactly one of the following holds.

\be\item $\emptyset\neq X^g$ --- that is, $g$ is elliptic.
\item
$\emptyset= X^g=A_g\subset\tilde{A_g}$, and $g$ is a nesting
reflection or an inversion. Moreover $\b_x(g)<0$ for all $x$ in this
case.
\item
$\emptyset=X^g\subset A_g=\tilde{A_g}$, and $g$ is hyperbolic.
Moreover $\b_x(g)>0$ for all $x$ in this case. \ee

\end{proposition}

\pf It is clear from the respective set inclusions that the cases
are mutually exclusive, and that $X^g\neq\emptyset$ precisely when
$g$ is elliptic. Moreover $A_g=\tilde{A_g}$ if $g$ is hyperbolic by
Lemma~\ref{nesting-lemma}(2) and proper inclusion holds if $g$ is a
nesting reflection or an inversion.

Fix $x\in X$, and put $u=\Y(g^{-1}x,x,gx)$. It is clear from
Theorem~\ref{hyper}(12) that $\b_x(g)=d(u,gu)>0$ for all $x$ if $g$
is hyperbolic. So suppose that $g$ is a nesting reflection or an
inversion. In this case $gu,g^{-1}u\in [x,u]$, for otherwise
(swapping $g$ and $g^{-1}$ if necessary) $gu\in[u,gx]$ which implies
that $[g^{-1}u;u;gu]$ which is impossible for a nesting reflection
or an inversion, by Lemma~\ref{nesting-lemma}(4). If
$[x,g^{-1}u,gu,u]$ then one can show that the following
configuration holds: $[g^{-1}x,g^{-2}u,u,g^2u,gu,g^{-1}u,x]$, and
hence $[g^{-1}x,u,gu,g^{-1}u,x]$. Using this and noting that
$gu=\Y(x,gx,g^{2}x)$, the defining expression for $\b_x(g)$
simplifies to $-d(u,gu)$. If instead $[x,gu,g^{-1}u,u]$, then the
foregoing shows that $\b_x(g^{-1})=-d(u,g^{-1}u)$. It can be shown
directly from the definition of $\b_x$ and using (L3) that
$\b_x(g)=\alpha_g\b_x(g^{-1})=-\alpha_g d(u,g^{-1}u)=-d(u,gu)$, as
required. \qed\\

Of course we can further distinguish between nesting reflections and
inversions: in the former case $X^{g^2}=\emptyset$ and $\b_x(g^2)>0$
for all $x$, while in the latter case $X^{g^2}\neq\emptyset$ and
$\b_x(g^2)=0$ for some $x$.

On the other hand, elliptic automorphisms cannot be distinguished
from other types of automorphism by considering $\b$ and $A_g$ and
$\tilde{A}_g$ alone.\\

Recall that if $g$ is a nesting reflection then $g^2$ is hyperbolic.
We thus have a useful criterion for recognising free actions.

\begin{corollary}\label{free-criterion}
\be
\item An affine action of a group on a $\Lambda$-tree is free and
without inversions if and only if $g^2$ is hyperbolic for all $g\neq
1$; \item If an affine action is free (possibly with inversions)
then $g^2$ is hyperbolic for all $g^2\neq 1$.
\item If an affine action is free, rigid and without inversions, then $g$ is hyperbolic for all $g\neq 1$.\qed\ee
\end{corollary}

To show that a free product of ATF groups is again ATF, we will
embed a family of ordered abelian groups in (a subgroup $\Lambda$
of) their Cartesian product and then use the base change functor to
obtain actions of the given groups on $\Lambda$-trees. The next step
is to show that freeness is preserved by the application of the base
change functor. However it is not true in general that the base
change functor preserves freeness, even if the map $h$ is an
embedding --- Example~\ref{nesting-ex} shows how it may fail (take
$\Lambda_1=\Lambda_0$ and $\Lambda_2=\R$, and take $h$ to be the
inclusion). However, in the case of interest to us this argument is
valid. To show this, we need a somewhat technical result which may
be of independent use.

Fix an $\alpha$-affine action of $G$ on a $\Lambda$-tree $X$, a
basepoint $x\in X$ and $g\in G$, and suppose that $[g^{-1}u;u;gu]$
where as usual $u=\Y(g^{-1}x;x;gx)$ (so $g$ is hyperbolic on the
subtree spanned by $\langle g\rangle x$). Note that $u\in A_g$. We
define the \emph{right radius of $A_g$ with respect to $x$} to be
$$\Rad^r(g)=\Rad_x^r(g)=\{d(u,v):v\in A_g,\ u\notin[gu,v]\}\cup\{0\}$$
The \emph{left radius of $A_g$ with respect to $x$}, $\Rad^l(g)$, is
the right radius of $A_{g^{-1}}$ with respect to $x$. Of course this
amounts to an abuse of notation, since these sets depend not just on
the set $A_g=A_{g^{-1}}$, but on $g$. However it is convenient in
this context to think of $A_g$ as having an orientation determined
by $g$ and we will do so. (Our terminology implies that $g$
translates `to the right', at least on the subtree spanned by
$\langle g\rangle u$.)

The \emph{diameter of $A_g$ with respect to $x$} is
$\Diam(g)=(-\Rad^l(g))\cup(\Rad^r(g))$. Clearly
$\Diam(g^{-1})=-\Diam(g)$. Moreover, $\Diam(g)$ is naturally
isometric to $A_g$, via a map which sends 0 to $u$.

In Example~\ref{nesting-ex}, if we take the subtree
$X_0=(-\infty,4/3]$ of $\R$, then $X_0$ is invariant under $g^2$.
Taking $x=u=1$, we see that $\Rad^r(g^2)=[0,\infty)$, and
$\Rad^l(g^2)=[0,1/3]$, so that $\Diam(g)=[-1/3,\infty)$. If instead
we define $X_0$ to exclude the fixed point $4/3$ then
$\Rad^l(g)=[0,1/3)$, so that $\Diam(g)=(-1/3,\infty)$.

\begin{proposition}\label{hyper-b}
Suppose that $X$ is spanned by the orbit $Gx$ and that $\b(g)> 0$
(where $\b=\b_x$). Then $g$ has a fixed point if and only if
$(1-\alpha_g)^{-1}\b(g)\cap\Rad^r(g)\neq\emptyset$ or
$(1-\alpha_{g^{-1}})^{-1}\b({g^{-1}})\cap\Rad^r({g^{-1}})\neq\emptyset$.
\end{proposition}

\pf Note first that $\b_x(g)=d(u,gu)$ if $\b_x(g)>0$; this can be
seen for example by applying Theorem~\ref{hyper}(12) to the subtree
of $X$ spanned by $\langle g\rangle x$.

Suppose that $y$ is fixed by $g$. Then $y\in A_g$ and, since
$\b(g)>0$, $y$ cannot lie in the subtree $T$ spanned by $\langle
g\rangle u$. Replacing $g$ by its inverse if necessary we have
$[u,gu,y]$, whence
\begin{eqnarray*} \alpha_g d(u,y) &=& d(gu,gy)\\ &=& d(gu,y)\\ &=&
d(u,y)-d(u,gu)\\ &=& d(u,y)-\b(g),
\end{eqnarray*}
forcing $d(u,y)\in(1-\alpha_g)^{-1}\b(g)$ and clearly
$d(u,y)\in\Rad^r(g)$.

Conversely, we have $[g^{-1}u;u;gu]$, so that
$\cup_{n\in\Z}[g^{n-1}u,g^nu]$ is a linear $\langle
g\rangle$-invariant subtree on which $g$ has no fixed point.
Replacing $g$ by its inverse if necessary, our hypothesis guarantees
the existence of $\lambda$ such that
$(1-\alpha_g)(\lambda)=\b(g)$ and  $\lambda$ lies in $\Rad^r(g)$. Since, moreover, $X$ is spanned by $Gx$,
$[g^{-1}u;u,y,\gamma x]$ for some $\gamma\in G$. We claim that
$T\cup\{g^{-1}y,y,gy\}$ is a collinear set. Since $T$ is $\langle
g\rangle$-invariant, it suffices to show that $T\cup\{g^{-1}y,y\}$
or $T\cup\{y,gy\}$ is collinear. If $y\in T$, the claim is obvious,
and if $y$ is in the subtree spanned by $\langle g\rangle x$ but
$y\notin T$ the configuration $[g^{-1}y,y,gy]$ is impossible. So
suppose otherwise. Since we have $[g^{\pm 1}t;t,y]$ for all $t\in
T$, there can be no point $x_1\in T$ such that $x_1\in[t,y]$ for all
$t\in T$. Therefore, applying Lemma~\ref{coll} to
$I_1=[g^{-1}y,y,gy]$ and $I_2=T$, there exists $x_0\in
[g^{-1}y,y,gy]$ such that $[t;x_0,g^{-1}y]$ and $[t;x_0,y]$ for all
$t\in T$. Thus either $\{t;x_0,g^{-1}y,y\}$ is collinear for all
$t\in T$ $x_0\in[g^{-1}y,y]$ so that $[g^{-1}y,x_0,y,gy]$ and
$[t,x_0,y,gy]$ is collinear for all $t\in T$. The claim follows.

We therefore have either $[u,y,gy]$ or $[u,gy,y]$, the configuration
$[y,u,gy]$ being forbidden by our choice of $y$. Now $gu\in[u,gy]$
and $\lambda=d(u,y)$ imply \begin{eqnarray*}\alpha_g(\lambda) &=&
d(gu,gy)\\ &=& d(u,gy)- d(u,gu).
\end{eqnarray*}

Therefore our choice of $\lambda$ gives
\begin{eqnarray*}\lambda=\b(g)+\alpha_g(\lambda) &=& d(u,gy)\\ &=&
d(u,y)\pm d(y,gy)\\ &=& \lambda\pm d(y,gy)
\end{eqnarray*}

This forces $y=gy$. \qed\\

The latter sentence of Proposition~\ref{hyper-b} may be restated as
follows: $g$ fixes a point if and only if
$(1-\alpha_g)^{-1}\b(g)\cap\Diam(g)\neq\emptyset$.

We conclude this discussion by observing that if $\b_x(g')>0$ for
all conjugates $g'$ of $g$ then $\Rad^r(g)$ (and hence $\Diam(g)$)
can be characterised in terms of based length functions $L=L_x$ in
the case where $X$ is spanned by $Gx$. Note first that $\b_{\gamma
x}(g)=\alpha_\gamma\b_x(\gamma^{-1}g\gamma)$ for $\gamma\in G$, so
that $\b_{\gamma x}(g)>0$ for all $\gamma$. So in the notation of
Theorem~\ref{hyper} if $w=u(\gamma x,g)$ then we have
$[g^{-1}w;w;gw]$. It follows that $\alpha_\gamma
\a_x(\gamma^{-1}g\gamma)=\a_{\gamma x}(g)=d(\gamma x,w)=d(\gamma
x,A_g)$, so that $\Rad^r(g)$ is spanned by elements of the form
$d(u,w)$ where $[x,u,w,\gamma x]$ and $gu$ and $\gamma x$ represent
the same direction at $u$ --- in other words
$\c_x(g,\gamma)>a_x(g)$. Thus

\begin{proposition}\label{radius}

$\Rad^r(g)$ is spanned by elements of the form
$$L_x(\gamma)-\a_x(g)-\alpha_\gamma \a_x(\gamma^{-1}g\gamma)\mbox{ where }\c_x(g,\gamma)>a_x(g).$$\qed\end{proposition}

\subsection{Line stabilisers}

A \emph{line} in a $\Lambda$-tree is a maximal linear subtree. Let
$G$ be a group that admits a free affine action (without inversions)
on a $\Lambda$-tree $X$. We will call a subgroup $H$ a \emph{line
stabiliser} for the action of $G$ if there is a line $T$ of $X$
whose stabiliser in $G$ is $H$.

Recall that the full group of isometric automorphisms of $\Lambda$
considered as a $\Lambda$-tree is isomorphic to $\Lambda\rtimes
C_2$: here $\Lambda$ acts on itself by translation, and the
non-trivial element of $C_2$ changes the sign (see \cite[2.5]{AB}).
It is not hard to generalise this to the affine case: the group of
affine automorphisms of $\Lambda$ itself is
$$\Aff(\Lambda)\cong\Lambda\ltimes\left(\Aut^+(\Lambda)\times
C_2\right).$$

(It is worth noting first that the set of permutations $g$ of a
$\Lambda$-tree $X$ that are $\alpha_g$-affine automorphisms for some
$\alpha_g\in\Aut^+(\Lambda)$ forms a group.)

Let $\Lambda$ be an ordered abelian group of finite rank. As noted
in \S1.2, if $\Lambda_0$ is a convex subgroup, order-preserving
automorphisms of $\Lambda_0$ extend to $\Lambda$, and
$\Aut^+(\Lambda_0)$ embeds in  $\Aut^+(\Lambda)$. Let $T$ be a
linear $\Lambda$-tree, which we assume to be isometrically embedded
in $\Lambda$, let $\alpha_g\in\Aut^+(\Lambda_0)$ where $\Lambda_0$
is spanned by $\{d(x,y):x,y\in T\}$, and $g$ an affine automorphism
of $T$ with dilation factor $\alpha_g$. Let
$\bar{\alpha}_g\in\Aut^+(\Lambda)$ be an extension of $\alpha_g$. We
claim that $g$ extends to an $\bar{\alpha}_g$-affine automorphism
$\Lambda$. First we can define an isometry $\sigma_g:T\to\Lambda$
via $x\mapsto(\bar{\alpha}_g)^{-1}gx$ It is easy to see that
$\sigma_g$ is isometric, and thus by \cite[Lemma
2.3.1]{Chiswell-book} $\sigma_g$ extends to an isometric
automorphism $\bar{\sigma}_g$ of $\Lambda$. It follows that
$\bar{\alpha}_g\bar{\sigma}_g$ is an $\bar{\alpha}_g$-affine
automorphism of $\Lambda$; it is routine to check that
$\bar{\alpha}_g\bar{\sigma}_g$ extends the action of $g$ on $T$.

We deduce that an affine action of a group on a linear subtree of $\Lambda$ extends to an affine action on $\Lambda$ provided $\Lambda$ has finite rank. (In fact the assumption of finite rank is only needed to ensure that the dilation factors extend to $\Lambda$.)

The following observation will also be useful. Let $G$ be a group
acting on a set $X$, $N$ a normal subgroup, and let $X_0=\{x\in
X:gx=x\ \forall g\in N\}$. Then $G$ stabilises $X_0$.

\begin{proposition}\label{sol}
Let $G$ be a group that admits a free affine action on a
$\Lambda$-tree $X$, and $H$ a subgroup of $G$.
\begin{enumerate}

\item
If $H$ is a non-trivial subnormal abelian subgroup of $G$ then $G$
stabilises a line.
\item
Suppose that $\Aut^+(\Lambda)$ is soluble and $H$ stabilises a line
of $X$. Then $H$ is soluble.
\item Suppose that $\Lambda=\Z^n$ for
some $n$ and $H$ stabilises a line $T$ of $X$ and preserves the
orientation. Then $H$ is nilpotent.
\end{enumerate}
\end{proposition}

\pf (1) For $g\neq 1$ we have $g^2$ hyperbolic by
Corollary~\ref{free-criterion}(1). Thus $g^2$ stabilises a line
($A_{g^2}$) and preserves the orientation of the line. If $\epsilon$
and $\epsilon'$ are the ends of this line then
$g^2\epsilon=\epsilon$ and $g^2\epsilon'=\epsilon'$, whence $g$
stabilises $\{\epsilon,\epsilon'\}$. Thus $g$ stabilises a line,
which must be unique since otherwise $g$ would fix a branch point.

If $gh=hg$ where $g$ and $h$ are non-trivial, then $g$ and $h$ must
stabilise a common line. It follows that a non-trivial abelian
subgroup $A$ of $G$ stabilises a unique line. The result follows by
the observation above and an easy induction on $n$.

(2) The subgroup $H$ acts faithfully on any line that it
stabilises, whence $H$ embeds in
$\Lambda\rtimes(\Aut^+(\Lambda)\times C_2)$. The result follows
since $\Lambda$ is abelian and $\Aut^+(\Lambda)$ is soluble.

(3) The full group of orientation-preserving affine automorphisms of
$T$ embeds in that of $\Z^n$, which is isomorphic to
$\Z^n\rtimes\Aut^+(\Z^n)\cong\UT(n+1,\Z)$, which is nilpotent. Thus,
since $H$ acts faithfully on $T$, $H$ is nilpotent. \qed

\begin{corollary}\label{linestabs}
Suppose that $G$ is non-trivial and has a free affine action on a
$\Lambda$-tree. \be\item If $\Aut^+(\Lambda)$ is soluble, the
non-trivial line stabilisers are precisely the maximal soluble
subgroups of $G$.
\item
If $\Lambda=\Z^n$ and there is no line whose ends are interchanged
by any $g\in G$, then the non-trivial line stabilisers are precisely
the maximal nilpotent subgroups of $G$.\ee
\end{corollary}

\pf (1) If $H\neq 1$ stabilises a line, then $H$ is soluble by
Proposition~\ref{sol}. If $\langle H,\gamma\rangle$ is soluble, then
by Proposition~\ref{sol} this subgroup stabilises a line, which must
be the same line as that stabilised by $H$. It follows that if $H$
is a line stabiliser then $H$ is maximal soluble.

Conversely, if $H$ is non-trivial and soluble, $H$ normalises the
last non-trivial term of its derived series, which is a non-trivial
abelian group, so that $H$ stabilises a line, by
Proposition~\ref{sol}. If $\gamma$ also stabilises this line, then
$\langle H,\gamma\rangle$ is soluble. Therefore if $H$ is maximal
soluble, then $H$ is a line stabiliser.

An analogous argument establishes (2). \qed\\

Note that the class of $\Lambda$ for which $\Aut^+(\Lambda)$ is
soluble includes all ordered abelian groups of finite rank: this
follows from Corollary~\ref{sol-f-rank}.

Recall that in the isometric case, if a group acts freely on a
$\Lambda$-tree $X$, and an end is fixed, then the group acts by
translations on an invariant subtree.
Thus each end stabiliser is a line stabiliser in the isometric case.
We next show that the same is true of affine actions on
$\Lambda$-trees provided $\Aut^+(\Lambda)$ is soluble. (In fact, the
statement is true under the slightly weaker hypothesis that the
image of $\alpha$ is soluble.)

\begin{proposition}
Let $\Lambda$ be an ordered abelian group\ for which
$\Aut^+(\Lambda)$ is soluble. Let $G$ be a group that has a free
affine action on a $\Lambda$-tree $X$, and suppose that $G$ fixes an
end $\epsilon$ of $X$. There is a unique line $(\epsilon,\epsilon')$
stabilised by $G$.
\end{proposition}

\pf Use induction on the derived length $d=d_G$ of
$\alpha(G)\leq\Aut^+(\Lambda)$. If $d=0$ then $G=\ker\alpha$ and the
given action is isometric. The proposition is true in the isometric
case, as we have already observed. Now if $[G,G]=1$ then $G$ is
abelian, and by Proposition~\ref{sol} $G$ stabilises a unique line.

Assume that the result true for groups $G_1$ with $d_{G_1}<d$, and
that $d\geq 2$. Then $[G,G]$ stabilises a unique line of $X$. Thus
$[G,G]$ fixes exactly two ends, one of which must be $\epsilon$, and
the other which we denote by $\epsilon'$.
Therefore $G$ stabilises this pair of ends; since $G$ fixes
$\epsilon$, it must also fix $\epsilon'$, whence $G$ stabilises the
line $(\epsilon,\epsilon')$.
Since $[G,G]$ fixes no other line, neither does $G$. \qed\\

Recall that a CSA group is one in which all
maximal abelian subgroups are malnormal. This means that for CSA
groups $G$, and maximal abelian subgroups $M$, if $M\cap
M^\gamma\neq 1$, then $\gamma\in M$. Note that this implies that
each non-trivial element of $G$ is contained in a unique
maximal abelian subgroup --- equivalently $G$ is commutative-transitive
(commutativity is an equivalence relation on its non-trivial
elements).

Let us consider these properties in a more general setting. Let
$\mathcal{P}$ be a class of groups that is subgroup-closed. Say that
a group $G$ is \emph{$\mathcal{P}$-disconnected} if each non-trivial
$x\in G$ is contained in a unique maximal $\mathcal{P}$
subgroup of $G$, called \emph{the $\mathcal{P}$-component of $x$}. (So abelian-disconnected groups are precisely
commutative-transitive groups.) Call a group $G$
$\mathrm{CS}\mathcal{P}$ if
$G$ is $\mathcal{P}$-disconnected and the maximal $\mathcal{P}$
subgroups of $G$ are malnormal.

Some cases of interest are where $\mathcal{P}$ is $A$,
the class of abelian groups, $S$, the class of soluble groups, and
$N$, the class of nilpotent groups.

\begin{proposition}\label{ATF-CSS}
\be\item Let $\Lambda$ be an ordered abelian group\ with
$\Aut^+(\Lambda)$ soluble. Then an $\mathrm{ATF}(\Lambda)$ group is
CSS.
\item
If $G$ admits a free affine action on a $\Z^n$-tree $X$ and no line
of $X$ has its ends interchanged by an element of $G$ then $G$ is
CSN.

\ee
\end{proposition}

\pf (1) and (2) may be proven by taking $\mathcal{P}=S$ and
$\mathcal{P}=N$, and invoking Corollary~\ref{linestabs}(1) and (2)
respectively. Let $M$ be a maximal $\mathcal{P}$ subgroup of a group
$G$ equipped with a free $\alpha$-affine action on a $\Lambda$-tree.
There exists a unique line $L=(\epsilon,\epsilon')$ in $X$ whose
stabiliser is $M$ by Corollary~\ref{linestabs}. Further,
the subgroup $\gamma M\gamma^{-1}$ is also maximal $\mathcal{P}$ and
stabilises $\gamma L$. If $h$ belongs to the intersection, then $h$
stabilises $L$ and $\gamma L$; but this forces $L=\gamma L$ since
$h$ stabilises a unique line. Thus $\gamma$ belongs to the
stabiliser of $L$, namely
$M$.\qed\\

We call a group $G$ an $\ATF[\mathcal{P}]$ group if it admits a free
affine action on a $\Lambda$-tree $X$ for some $\Lambda$ such that
all line stabilisers are $\mathcal{P}$ subgroups, and every
$\mathcal{P}$ subgroup of $G$ stabilises a line (which must be unique). The argument of Proposition~\ref{ATF-CSS} shows that $\ATF[\mathcal{P}]$ groups are CS$\mathcal{P}$ (and hence $\mathcal{P}$-disconnected). The following
proposition follows the idea of \cite[Theorem 5.5.3]{Chiswell-book}
(which is largely due to B. Baumslag and Remeslennikov) in the case
where $\mathcal{P}$ is the variety of abelian groups, and free
groups are considered in place of $\ATF[\mathcal{P}]$ groups. In
this case the assumption that $G$ is not a $\mathcal{P}$ group can
be dropped.

The following proposition is in fact true if $\ATF[\mathcal{P}]$ is replaced by any group theoretic property stronger than CS$\mathcal{P}$.

\begin{proposition}
Let $\mathcal{P}$ be a variety of groups that contains all abelian groups, and suppose that $G$ is residually $\ATF[\mathcal{P}]$, but that $G$ is not a $\mathcal{P}$ group. The following are equivalent.
\be
\item $G$ is fully residually $\ATF[\mathcal{P}]$.
\item $G$ is 2-residually $\ATF[\mathcal{P}]$.
\item $G$ is $\mathcal{P}$-disconnected.
\item $G$ is CS$\mathcal{P}$.
\ee
\end{proposition}

\pf
(1)$\Rightarrow$(2) is trivial. Suppose that (2) holds. If $H$ is a $\mathcal{P}$ subgroup of $G$ then any homomorphism $G\to Q$ maps $H$ to a $\mathcal{P}$ subgroup of $Q$. Let $x,y,z\in G$ with $y\neq 1$ and suppose that $\langle x,y\rangle$ and $\langle y,z\rangle$ are $\mathcal{P}$ subgroups: we claim that $H=\langle x,y,z\rangle$ is $\mathcal{P}$. For taking any law $w=w(x_1,\ldots,x_r)$ of $\mathcal{P}$, let $g=w(g_1,\ldots, g_r)$ where $g_i\in H$ for all $i$. If $g\neq 1$ we can apply (2) to obtain a homomorphism $\phi:G\to Q$ where $Q$ is $\ATF[\mathcal{P}]$, and hence $\mathcal{P}$-disconnected, and where $\phi(g)\neq 1\neq\phi(y)$. Now $\phi(x)$ and $\phi(z)$ both belong to the maximal $\mathcal{P}$ subgroup containing $\phi(y)$. Hence $\phi(g)=w(\phi(g_1),\ldots,\phi(g_r))=1$, a contradiction. Thus $g=1$ from which we conclude that $H$ is a $\mathcal{P}$ subgroup.

Next let $M=M_y=\{x\in G:\langle x,y\rangle\mbox{\ is a }\mathcal{P}\mbox{ subgroup}\}$. It is now straightforward to show that $M$ is in fact a subgroup. Moreover, every pair of (non-trivial) elements of $M$ is mapped to the same $\mathcal{P}$ component by any homomorphism with $\mathcal{P}$-disconnected codomain, which forces $M$ to be a $\mathcal{P}$ subgroup, using the 2-residual $\mathcal{P}$-disconnectedness of $G$. Clearly $\langle M,\gamma\rangle$ can only be a $\mathcal{P}$ subgroup if $\gamma\in M$, and any $\mathcal{P}$ subgroup of $G$ containing $y$ is contained in $M$; therefore $M$ is the unique maximal $\mathcal{P}$ subgroup of $G$ containing $y$. This proves (3).

Assume (3). Let $M$ be a maximal $\mathcal{P}$ subgroup of $G$, and suppose that $\gamma\notin M$. Then $\langle M,\gamma\rangle$ is not a $\mathcal{P}$ subgroup and there exists a law $w$ for the variety $\mathcal{P}$ such that $g=w(g_1,\ldots,g_r)\neq 1$ for some choice of $g_i\in \langle M,g\rangle$ ($1\leq i\leq r$). Let $Q$ be an $\ATF[\mathcal{P}]$ group and $\phi:G\to Q$ a homomorphism with $\phi(g)\neq 1$. Now $\phi(M)$ is a $\mathcal{P}$ subgroup of $Q$ which is thus contained in a unique maximal $\mathcal{P}$ subgroup $\bar{M}$, since $\ATF[\mathcal{P}]$ groups are $\mathcal{P}$-disconnected.
If $\phi(\gamma)\in\bar{M}$, then $\langle\phi(M),\phi(\gamma)\rangle$ is a $\mathcal{P}$ subgroup of $Q$. But then $\phi(g)=w(\phi(g_1),\ldots,\phi(g_r))$ with each $\phi(g_i)$ belonging to the $\mathcal{P}$ subgroup $\langle\phi(M),\phi(\gamma)\rangle$. This forces $\phi(g)=1$, a contradiction.

Therefore $\phi(\gamma)\notin\bar{M}$, giving $\bar{M}\cap\phi(\gamma)\bar{M}\phi(\gamma^{-1})=1$ since $Q$ is CS$\mathcal{P}$. Suppose that $M$ and $\gamma M\gamma^{-1}$ have non-trivial intersection. Since these are maximal $\mathcal{P}$ subgroups and $G$ is $\mathcal{P}$-disconnected, they must coincide. But $\phi(M)\subseteq\bar{M}$ and $\phi\left(\gamma M\gamma^{-1}\right)\subseteq\phi(\gamma)\bar{M}\phi(\gamma^{-1})$, which gives a contradiction.
Therefore $M\cap \gamma M\gamma^{-1}=1$, whence (4).

Finally assume (4). Then for non-trivial $x\in G$, we can find $g\in G$ such that $g$ does not belong to the $\mathcal{P}$ component $M$ of $x$. Thus $M\cap gMg^{-1}=1$. Let $y$ be another non-trivial element of $G$. If $[x,y]=1$ then $x$ and $y$ must belong to the same $\mathcal{P}$ component of $G$. We therefore have $gxg^{-1}$ and $y$ in distinct $\mathcal{P}$ components.
In particular $[gxg^{-1},y]\neq 1$.

Now let $x_1,\ldots x_n$ be given non-trivial elements of $G$. Choosing a suitable conjugate $y_i$ of $x_i$ as above we have $[y_1,y_2,\ldots,y_n]\neq 1$. Using the residual $\ATF[\mathcal{P}]$ property we can map $G$ to an $\ATF[\mathcal{P}]$ group $Q$ with the image of $[y_1,y_2,\ldots,y_n]$ non-trivial. This forces the image of each $y_i$ to be non-trivial, and with it, the image of each $x_i$. This proves (1).
\qed

\section{Examples}

Of course all isometric actions are affine --- this corresponds to
the case where the given homomorphism $\alpha$ is trivial.
Consequently, $\ITF(\Lambda)$ implies $\ATF(\Lambda)$ as already
noted.

The group $\Hol^+(\Lambda)=\Lambda\rtimes\Aut^+(\Lambda)$ acts
naturally on $\Lambda$: this action is affine though not free in
general.

\example[Some free affine actions on $\R$-trees]\label{liousse-exx}
In \cite{Liousse}, I. Liousse constructs examples of groups that
admit free affine actions on $\R$-trees. She constructs two types,
the first having presentations of the form
$$\langle x_1,x_2,\ldots, x_n,y_1,y_2,\ldots,y_n\ |\
[x_1,y_1]=[x_2,y_2]=\cdots=[x_n,y_n]\rangle,$$

and the second $$\langle x_1,x_2,\ldots, x_n,y_1,y_2,\ldots,y_m\ |\
w=v^k\rangle$$ where $w$ is a surface relator (that is, either
$w=[x_1,x_2]\cdots[x_{2r-1},x_{2r}]$ where $n=2r$, or $w=x_1^2\cdots
x_n^2$) and $v$ is a non-trivial word in the $y_i^{\pm}$.

As Liousse notes, it can be shown using Rips' Theorem that most of these groups do not admit free isometric actions on $\R$-trees. This shows, in our notation, that $\ATF(\R)$ does not imply $\ITF(\R)$.\\

%
%
\example[A non-free affine action of a free group] Let $F$ be the
free group on $\{a_n:n\in\Z\}$, and $\Gamma$ the corresponding
Cayley graph. (So there is an edge labelled $a_n$ joining $\gamma$
to $\gamma a_n$ for $\gamma\in F$ and $n\in\Z$). Define the distance
between two adjacent vertices ($=$ length of the edge) to be $2^{n}$
if $a_n$ is the label of the edge joining them. Extend to an
$\R$-metric on the vertex set of $\Gamma$ by putting $d(x,y)$ equal
to the sum of the lengths of the edges in the reduced path joining
$x$ and $y$. This is clearly a 0-hyperbolic $\R$-metric.

Now left multiplication by $F$ amounts to an isometric action on
$\Gamma$. Moreover, we can define a `shifting' map $\tau:F\to F$ by
mapping an element $a_{i_1}^{\epsilon_1}a_{i_2}^{\epsilon_2}\cdots
a_{i_m}^{\epsilon_m}$ to
$a_{i_1+1}^{\epsilon_1}a_{i_2+1}^{\epsilon_2}\cdots
a_{i_m+1}^{\epsilon_m}$. It is routine to verify that $d(\tau x,\tau
y)=2d(x,y)$ for $x,y\in F$, whence $G=\langle F,\tau\rangle$ has an
affine action on $\Gamma$: the map $\alpha:G\to\Aut^+(\R)$ is given
by $\gamma\mapsto 2^{e}$ where $e$ is the exponent sum of $\tau$ in
$\gamma$.

Now $\Gamma$ embeds in an $\R$-tree in a natural way by
Theorem~\ref{basechange}. Moreover the action of $G$ on $\Gamma$
extends to an affine action on the $\R$-tree.

Observe next that
$\tau^{-1}a_k\tau(a_{i_1}^{\epsilon_1}a_{i_2}^{\epsilon_2}\cdots
a_{i_m}^{\epsilon_m})=a_{k-1}(a_{i_1}^{\epsilon_1}a_{i_2}^{\epsilon_2}\cdots
a_{i_m}^{\epsilon_m})$. Suppose that $w$ is a word in $\tau$ and
$a_k$ ($k\in\Z$) which acts as the identity. Then $w\in\ker\alpha$,
which means that the exponent sum of $\tau$ in $w$ is zero.
Therefore $w$ is expressible as a product of terms of the form
$\tau^{-k}a_l^{\pm 1}\tau^k$. Using the relation scheme noted above,
$w$ is in fact expressible as a product of $a_{l'}^{\pm 1}$
($l'\in\Z$); that is, $w\in F$. This forces $w=1$, since the action
of $F$ on its Cayley graph is free.

Therefore $G$ has the presentation $$\langle \tau, a_k\  (k\in\Z)\ |\ \tau^{-1}a_k\tau=a_{k-1}\rangle\cong\langle a_0,\tau\ |\ \rangle,$$ that is, $G$ is the free group on $\{a_0,\tau\}$.\\

%
%
\example[The Heisenberg group is $\ATF$] Let $X=\Lambda=\Z^3$, and
let $\sigma:(x,y,z)\mapsto(x,y+1,z)$,
$\tau:(x,y,z)\mapsto(x+1,y,z+y)$. Let $G=\langle \sigma,
\tau\rangle$. Then $G$ has an $\alpha$-affine action on $X$ where
$\alpha_\sigma=\id$ and $\alpha_\tau(x,y,z)=(x,y,z+y)$.

Suppose that $w\in G$ maps $[r,s]$ into itself for some $r\leq s\in X$, so that (replacing $w$ by $w^2$ if necessary) $r\leq wr\leq ws\leq s$. Then the exponent sum of $\tau$ in $w$ must be zero
so that $w$ is a product of conjugates of the form $\tau^{-k}\sigma^l\tau^k$, each of which fixes the first entry of each element of $X$, and adds $l$ to the second entry. Thus the exponent sum of $\sigma$ in $w$ is also zero, so that $w$ lies in the derived subgroup of $\langle \sigma,\tau\rangle$. Direct
computation shows that $[\sigma,[\sigma,\tau]]=[\tau,[\sigma,\tau]]=1$; it follows that $w$ is a product of commutators of the form $[\sigma^l,\tau^k]$. It is straightforward to show that $[\sigma^l,\tau^k](x,y,z)=(x,y,z-kl)$ and thus $w$ is a product of commutators of the form $[\sigma^{l_i},\tau^{k_i}]$ with $\sum k_il_i=0$, which therefore fixes $X$ pointwise. This shows that the action is free and rigid.

The calculation above also implies
$$[\sigma,\tau]^{kl}=[\sigma^l,\tau^k],\h k,l\in\Z.$$
It follows that $$G\cong\langle \sigma, \tau\ |\
[\sigma,\tau]^{kl}=[\sigma^l,\tau^k],\ k,l\in\Z\rangle.$$ We claim
that in fact $G$ is isomorphic to $H=H_3(\Z)$, the discrete
Heisenberg group. It is well known that $H$ is generated by the
matrices $x=\left(\begin{array}{ccc}1 & 1 & 0 \\ 0 & 1 & 0 \\ 0 & 0
& 1\end{array}\right)$ and $y=\left(\begin{array}{ccc}1 & 0 & 0
\\ 0 & 1 & 1 \\ 0 & 0 & 1\end{array}\right)$, and that
$[x,[x,y]]=[y,[x,y]]=1$ are defining relations.

Since $[\sigma,[\sigma,\tau]]=[\tau,[\sigma,\tau]]=1$, there
is an epimorphism $\phi:H\to G$ defined by $x\mapsto\sigma$ and
$y\mapsto\tau$. Moreover,
$[x^k,y^l]=[x,y]^{kl}$ for all $k,l\in\Z$, whence $\sigma\mapsto x$
and $\tau\mapsto y$ defines an epimorphism $\psi:G\to H$. Clearly
$\phi\circ\psi=\id$, whence both maps are isomorphisms.

Thus $H$ admits a free rigid affine action on $\Z^3$, viewed as a
$\Z^3$-tree, as claimed in Theorem~\ref{main-exx}(1).
\\

%
%
\example[$C_{\infty}\wr C_{\infty}$ and the soluble Baumslag-Solitar
groups are $\ATF$]

Let $X=\Lambda=\Z\times\R$, and $a$ a positive real number. Define
$\sigma=\sigma_a:(x,y)\mapsto(x+1,ay)$ and
$\tau:(x,y)\mapsto(x,y+1)$ and let $G_a=\langle \sigma_a,
\tau\rangle$. Then
$\sigma^{-k}\tau^l\sigma^k(x,y)=(x,y+\frac{l}{a^k})$. It follows
that $$[\sigma^{-k}\tau^l\sigma^k,
\sigma^{-k'}\tau^{l'}\sigma^{k'}]=1.$$

Let us write $C$ for the set consisting of these relations (where
$k,k',l,l'$ range through the integers). Taking $H=\langle
\tau\rangle$, and $K=\langle \sigma\rangle$ acting on the product
$\prod_{n\in\Z}\sigma^{-n} H\sigma^{n}$ by conjugation, we see that
$G_a$ is a quotient of $H\wr K\cong \langle \sigma, \tau|\
C\rangle\cong C_\infty\wr C_\infty$.

Suppose that $w$ is a word in $\sigma$ and $\tau$ and that
$w[r,s]\subseteq[r,s]$. Then the exponent sum of $\sigma$ in $w$
must be zero, so that $w$ is a product of conjugates of $\tau$ by
powers of $\sigma$. Since  these conjugates commute as observed
above we can write
$w=\sigma^{-k_1}\tau^{l_1}\sigma^{k_1}\cdots\sigma^{-k_r}\tau^{l_r}\sigma^{k_r}$
where $k_i< k_{i+1}$ for each $i$. Since the effect of such a $w$ is
to add $\sum_{i=1}^r l_ia^{-k_i}$ to the second entry, we must have
$\sum_{i=1}^r l_ia^{-k_i}=0$. This also shows that the action of
$G_a$ is free. Multiplying both sides of this equation by $a^{k_r}$,
we obtain a polynomial equation. We next consider the consequences
of this equation. Distinguish three cases.

\begin{itemize}
\item If $a$ is transcendental, then all $l_i$ are equal to zero since otherwise the polynomial equation above has $a$ as a root. Therefore
    $C$ amounts to a set of defining relations for $G_a$ in this case, so that $G_a\cong C_\infty\wr C_\infty$.

\item
    For any polynomial $p(x)\in\Z[x]$ for which $p(a)=0$, we obtain a relation for the group $G_a$, as follows. Write $p(x)=\sum_{k=0}^nl_kx^k$ and put $$w_p=
\tau^{l_0}\cdot(\sigma\tau^{l_1}\sigma^{-1})\cdots(\sigma^{n}\tau^{l_n}\sigma^{-n}).$$
Then $w_p=1$ is a relation of $G_a$. It is easily checked that if
$p_1(x)=dx^m$, then
$w_{p_1p}=\sigma^m\left(w_p\right)^d\sigma^{-m}$, which is clearly
equal to the identity in $G_a$ if $w_p=1$. Moreover, if
$p_1(x),p_2(x)\in\Z[x]$ then $w_{(p_1+p_2)}$ is easily seen to
 be equal to $w_{p_1}w_{p_2}$. It follows that if $p_0(a)=0$ and $p_0$ is a divisor of $p$, then $w_p$ is a consequence of $w_{p_0}$ and $C$.

 Conversely, it is clear that if $w$ is a word in $\sigma$ and $\tau$ giving the identity in $G_a$, $w$ must have the form $w_p$ for some $p(x)\in\Z[x]$ with $p(a)=0$.
Therefore taking $p_0(x)$
 to be the minimum polynomial of $a$,
 we obtain the following presentation.
 $$G_a\cong\langle\sigma,\tau\ |\ w_{p_0}=1,\
 C \rangle$$

    Of course if $a_1$ and $a_2$ have the same minimum polynomial, then $G_{a_1}$ and $G_{a_2}$ are isomorphic.

\item If $a$ is a (positive) integer, then
$$G_a\cong\langle \sigma,\tau\ |\ \sigma \tau\sigma^{-1}=\tau^a,\ C\rangle.$$

In fact the first relation implies that
$\sigma^{k}\tau\sigma^{-k}=\tau^{a^{k}}$ and
$\sigma^{k'}\tau\sigma^{-k'}=\tau^{a^{k'}}$ commute for $k,k'\geq 0$
--- and hence for all $k,k'\in\Z$. It follows that $$G_a\cong \
\langle\sigma,\tau\ |\ \sigma \tau\sigma^{-1}=\tau^a\rangle \cong
\mathrm{BS}(1,a),$$ a soluble Baumslag-Solitar group.

If $a^{-1}$ is an integer, then replacing $\tau$ by its inverse we
also find that $G_a$ is isomorphic to $\mathrm{BS}(1,a)$.

\end{itemize}

The group $G_a$ for any transcendental $a$ gives the example
promised by Theorem~\ref{main-exx}(2), while taking $a\in\N$ gives
the example referred to in Theorem~\ref{main-exx}(3). Moreover all
actions under consideration here are in fact rigid.

\section{Constructions}

We say that a length function $L$ on a group $G$ is \emph{regular}
if for all $g,h\in G$, there exists $u\in G$ and $g',h'\in G$ with
$L(u)=\c(g,h)$ where $L(g)=L(u)+\alpha_uL(g')$ and
$L(h)=L(u)+\alpha_uL(h')$, and $g= ug'$, $h= uh'$. In the isometric
case, this boils down to the notion of regular length function in
the sense of Myasnikov and Remeslennikov. Call an action regular if
there exists a basepoint with respect to which the associated based
length function is regular.

\subsection{Free products}

If $\gamma$ is an element of a free product, we use the notation
$\gamma=\gamma_1\cdot\gamma_2\cdots\gamma_k$ if the last syllable of
each $\gamma_i$ (when $\gamma_i$ written as a reduced word in
elements of the free factors) belongs to a different free factor
from the first syllable of $\gamma_{i+1}$. (We do not assume that
each $\gamma_i$ belongs to a free factor.) Equivalently, the
syllable length of $\gamma$ as an element of the free product is
equal to the sum of the syllable lengths of the $\gamma_i$.
\\

\begin{theorem}\label{freeproduct} \emph{(see
Theorem~\ref{main-free-product})}\\ Let $G_i$ be a group admitting
an affine action on a $\Lambda$-tree for each $i\in I$. The free
product $G=\ast_{i\in I}G_i$ has an affine action on a
$\Lambda$-tree which extends the given actions of $G_i$.

If the given actions of $G_i$ are respectively \be\item free, \item
free and rigid,
\item or regular\ee then so is that of $G$.

\end{theorem}

\pf Suppose that for each $i$ we have an $\alpha_i$-affine action on
a $\Lambda$-tree. Take a basepoint $x_i$ from each $\Lambda$-tree,
and let $L=L_{x_i}$ denote the associated based length function.

Let $\alpha:G\to\Aut^+(\Lambda)$ be the unique common extension of
the $\alpha_i$. If $g_1,g_2,\ldots g_n$ are group elements, we
define $\overline{g_k}=g_1\cdots g_k$ for $1\leq k\leq n$ and
$g_0=\overline{g_0}=1$. It is also convenient to put $g_m=1$ for
$m>n$.

Now for $g=\overline{g_n}$ with consecutive $g_i$ in distinct free
factors, we put $$L(g)=\sum_{k=1}^n
\alpha_{\overline{g_{k-1}}}L(g_k).$$

(Here we abuse notation by using the same notation $L$ for the
length function on each $G_i$ and for the length function on $G=\ast
G_i$; similar abuses follow regarding the ancillary functions $\b$
and $\c$. In each case the function defined on $G$ restricts to the
given function on each $G_i$.)

It is clear that $L(1_G)=0$. Moreover,
$g^{-1}=g_n^{-1}g_{n-1}^{-1}\cdots g_1^{-1}$, and
$$\begin{array}{rl}\alpha_g L(g^{-1}) &=\alpha_g\sum_{k=1}^n \alpha_{g_n^{-1}\cdots g_{n-k+2}^{-1}} L(g_{n-k+1}^{-1})\\ &=
\sum_{k=1}^n \alpha_{\overline{g_{n-k}}}\
\alpha_{g_{n-k+1}}L(g_{n-k+1}^{-1})\\ & =\sum_{k=1}^n
\alpha_{\overline{g_{n-k}}}L(g_{n-k+1})\\ &
=\sum_{k=1}^n\alpha_{\overline{g_{k-1}}}L(g_k)\\ &=L(g)\end{array}$$

\begin{claim}
Let $\gamma=\gamma_1\cdot\gamma_2\cdots\gamma_q$. Then
$L(\gamma)=\sum_{k=1}^q
\alpha_{\overline{\gamma_{k-1}}}L(\gamma_k)$. In particular if
$L(\gamma)=0$ then $L(\gamma_i)=0$ for all $i$.
\end{claim}

The claim is trivial if $q=1$. If $q=2$, and $\gamma_1=g_1\cdots
g_p$ and $\gamma_2=g_{p+1}\cdots g_r$, with $g_i$ and $g_{i+1}$
belonging to distinct free factors for each $i$, it is
straightforward to show that both sides of the desired equation are
equal to $$\sum_{k=1}^r \alpha_{\overline{g_{k-1}}}L(g_k).$$ The
claim follows by an easy induction on $q$.

\begin{claim}
Suppose that $g=g_1\cdot g_2\cdots g_n$ and $h=h_1\cdot h_2\cdots
h_m$, and that $p\geq 0$ is the integer satisfying $g_i=h_i$ for
$1\leq i\leq p$, and either $g_{p+1}\neq h_{p+1}$ or
$p=\min\{m,n\}$. Put $g_i=1$ for $i>n$ and $h_i=1$ for $i>m$. Then
$\c(g,h)=L(\overline{g_p})+\alpha_{\overline{g_p}}\c(g_{p+1},h_{p+1})$.

If $g_{i}$ and $h_{i}$ belong to distinct free factors,
$\c(g_{i},h_{i})=0$.
\end{claim}

Note first that if $g_i$ and $h_i$ belong to distinct free factors,
we have \begin{eqnarray*}\alpha_{g_i}L(g_i^{-1}\cdot
h_i)&=&\alpha_{g_i}\left(L(g_i^{-1})+\alpha_{g_i^{-1}}L(h_i)\right)\\
&=&L(g_i)+L(h_i).\end{eqnarray*}

Thus $\c(g_i,h_i)=0$.

Write $g^{-1}h=g_n^{-1}\cdot g_{n-1}^{-1}\cdots
g_{p+2}^{-1}\cdot\left(g_{p+1}^{-1}h_{p+1}\right)\cdot h_{p+2}\cdots
h_m$.

Using Claim~34, it is then straightforward to show that
\begin{eqnarray*}\alpha_gL(g^{-1}h)&=&\alpha_{\overline{g_{p+1}}}\left(L(g_{p+2}\cdots
g_n)+L(g_{p+1}^{-1}h_{p+1})+
\alpha_{{g_{p+1}^{-1}}{h_{p+1}}}L(h_{p+2}\cdots
h_m)\right).\end{eqnarray*}

Now \begin{eqnarray*}L(g)&=&L(\overline{g_{p+1}})+\alpha_{\overline{g_{p}}g_{p+1}}L(g_{p+2}\cdots g_n) \\
\mbox{and\ \
}L(h)&=&L(\overline{h_{p+1}})+\alpha_{\overline{h_{p}}h_{p+1}}L(h_{p+2}\cdots
h_m).\end{eqnarray*} Combining these expressions, and recalling that
$\overline{g_p}=\overline{h_p}$, we obtain
\begin{eqnarray*}2\c(g,h)&=&L(\overline{g_{p+1}})+L(\overline{h_{p+1}})
-\alpha_{\overline{g_{p+1}}}L(g_{p+1}^{-1}h_{p+1})\\
&=&L(\overline{g_{p}})
+\alpha_{\overline{g_p}}L(g_{p+1})+L(\overline{h_p})+\alpha_{\overline{h_p}}L(h_{p+1})
-\alpha_{\overline{g_{p+1}}}L(g_{p+1}^{-1}h_{p+1})\\ &=&
2L(\overline{g_{p}})+2\alpha_{\overline{g_p}}\c(g_{p+1},h_{p+1})
\end{eqnarray*}
which establishes the claim.\\

It is now clear that $\c(g,h)\in\Lambda$ for all $g,h\in G$.

We now verify the isosceles condition for $\c(g,h)$, $\c(g,k)$ and
$\c(h,k)$. Keeping the notation for $g$ and $h$ as in Claim 35,
write $k=k_1\cdot k_2\cdots k_l$. Let $q\geq0$ be the integer for
which $h_i=k_i$ for $0\leq i\leq q$, and either $q=\min\{m,l\}$ or
$h_{q+1}\neq k_{q+1}$. Assume without loss of generality that $p\leq
q$. (This means that $h$ and $k$ share a common initial subword that
is no shorter than that of $g$ and $h$, when these elements of $G$
are written as words in the free factors.) Observe that $g_i=k_i$
for $i\leq p$ and $g_{p+1}\neq k_{p+1}$, and therefore that
$\overline{g_p}=\overline{h_p}=\overline{k_p}$.

Now Claim 35 gives \begin{eqnarray*}\c(g,h)&=&L(\overline{g_p})+\alpha_{\overline{g_p}}\c(g_{p+1},h_{p+1})\\ \c(g,k)&=&L(\overline{g_p})+\alpha_{\overline{g_p}}\c(g_{p+1},k_{p+1})\\
\mbox{and\ \ }
\c(h,k)&=&L(\overline{h_q})+\alpha_{\overline{h_q}}\c(h_{q+1},k_{q+1})
.\end{eqnarray*}

Suppose first that $\c(g,h)>\c(h,k)$. Then
$L(\overline{g_p})+\alpha_{\overline{g_p}}\c(g_{p+1},h_{p+1})
>L(\overline{g_p})+\alpha_{\overline{g_p}}L(h_{p+1}\cdots h_q)+\alpha_{\overline{h_q}}\c(h_{q+1},k_{q+1})$. If $q=p$ this inequality reduces to $\c(g_{p+1},h_{p+1})>\c(h_{p+1},k_{p+1})$, forcing $\c(h_{p+1},k_{p+1})=\c(g_{p+1},k_{p+1})$ and thus $\c(h,k)=\c(g,k)$.

Otherwise $q\geq p+1$, which implies
$\c(g_{p+1},h_{p+1})>L(h_{p+1}\cdots h_q)\geq L(h_{p+1})$,
contradicting Lemma~\ref{length-props}(8).

Next suppose that $\c(g,k)>\c(g,h)$. Then
$\c(g_{p+1},k_{p+1})>\c(g_{p+1},h_{p+1})$. Therefore the left-hand
side is non-zero, which forces $g_{p+1}$ and $k_{p+1}$ to belong to
the same free factor, by Claim 35. If $h_{p+1}$ belongs to a
different free factor, then $h_{p+1}\neq k_{p+1}$ giving $q=p$. Thus
$\c(h_{p+1},k_{p+1})=0=\c(g_{p+1},h_{p+1})$, which implies
$\c(h,k)=\c(g,h)$.

If $h_{p+1}$ belongs to the same free factor as $g_{p+1}$ and
$k_{p+1}$ then
$\c(h_{p+1},k_{p+1})=\c(g_{p+1},h_{p+1})<\c(g_{p+1},k_{p+1})$. Again
we must have $q=p$, for otherwise $q\geq p+1$, so that $g_{p+1}\neq
k_{p+1}=h_{p+1}$ and thus
$L(k_{p+1})\geq\c(g_{p+1},k_{p+1})>\c(h_{p+1},k_{p+1})=L(k_{p+1})$
by Lemma~\ref{length-props}(5), a contradiction. The equality
$\c(h_{p+1},k_{p+1})=\c(g_{p+1},h_{p+1})$ now gives the required
$\c(h,k)=\c(g,h)$.

Finally, suppose that $\c(h,k)>\c(g,h)$. If $q\geq p+1$ then
$h_{p+1}=k_{p+1}$
giving $\c(g,h)=\c(g,k)$. Otherwise $q=p$ and this inequality yields
$L(\overline{h_p})
+\alpha_{\overline{h_p}}\c(h_{p+1},k_{p+1})>L(\overline{g_p})+\alpha_{\overline{g_p}}\c(g_{p+1},h_{p+1})$,
giving
$\c(h_{p+1},k_{p+1})>\c(g_{p+1},h_{p+1})=\c(g_{p+1},k_{p+1})$. This
implies $\c(g,h)=\c(g,k)$.

The remaining cases, such as $\c(g,k)>\c(h,k)$, follow from those we
have considered by swapping the roles of $h$ and $k$.

This completes the proof that $L$ is an $\alpha$-affine length function.\\

We now prove that if the given actions of $G_i$ are free and without
inversions, then so is that of $G$. This follows easily from our
final claim:

\begin{claim}
Suppose that $r>1$, and $g=g_1\cdot\cdots\cdot g_r$ is cyclically
reduced (i.e. $g_1$ and $g_r$ belong to distinct free factors). If
$L(g_i)\neq 0$ for all $i$ then $g$ is hyperbolic and rigid.
\end{claim}

We have $L(g)\neq 0$. For otherwise, by Claim 34, if
$g=g_1\cdot\cdots\cdot g_r$ then $L(g_i)=0$ for all $i$,
contradicting the non-triviality of $g$. Let $x$ be the basepoint of
the $\Lambda$-tree $X$ corresponding to the length function $L$ on
$G$. Observe that if $g^nx=g^mx$ for some $m<n$ then $L(g^k)=0$ for
some $k\geq 1$, giving $0=L(g^k)=L(g)+\alpha_g
L(g)+\cdots+\alpha_{g^{n-1}}L(g)\geq L(g)$ --- here we have used the
assumption that $g$ is cyclically reduced and has length at least 2.
Thus $L(g)=0$, which is impossible. Note also that the equality
$$L(g^2)=L(g\cdot g)=L(g)+\alpha_g L(g),$$ which follows from Claim
1, implies $[x,gx,g^2x]$, giving $x\in A_g$. Since
$A_g\neq\emptyset$, $g$ cannot be a nesting reflection or an
inversion. This means that $g$ restricted to the invariant subtree
$\cup_{n\in\Z}[g^{n-1}x,g^nx]\subseteq A_g$ is hyperbolic. To
establish hyperbolicity of $g$ as an automorphism of $X$, it
suffices to show this set in fact coincides with $A_g$.

So suppose otherwise. Observe that $X$ is spanned by points of the
form $\gamma x$ where $\gamma$ ranges through $G$. Replacing $g$ by
its inverse if necessary, our supposition gives, for some $\gamma\in
G$,
$$[x,g^nx,\gamma x]\h\forall
n\geq 1.$$ Now
\begin{eqnarray*}[x,\gamma_1x,\gamma_2x]=[x,\gamma_2x]
&\Leftrightarrow&L(\gamma_1)+\alpha_{\gamma_1}L(\gamma_1^{-1}\gamma_2)=L(\gamma_2)\\
&\Leftrightarrow& L(\gamma_1)=\c(\gamma_1,\gamma_2)\end{eqnarray*}
for $\gamma_1,\gamma_2\in G$, so replacing $\gamma_1$ by $g^n$ and
$\gamma_2$ by $\gamma$ we get
\begin{equation}\label{cgnl}\c(g^n,\gamma)=L(g^n)\h\forall n\geq
1.\end{equation}

To show that $g$ is hyperbolic (with respect to the action on $X$),
it suffices to show that $g^n$ is hyperbolic for some $n$, since $g$
is hyperbolic with respect to the action on
$\cup_{n\in\Z}[g^{n-1}x,g^nx]$, and thus cannot act as an inversion
or a nesting reflection. So if $g$ is an initial subword of $\gamma$
(i.e. if $\gamma$ can be written in the form $g\cdot g^{-1}\gamma$),
then we can replace $g$ by a sufficiently large power of $g$ to
ensure that this is not the case. Now note that the longest common
initial subword $\overline{g_p}$ of $g$ and $\gamma$ is the same as
that $g^n$ and $\gamma$. So by Claim 35, we have
\begin{equation}\label{cgn}\c(g^n,\gamma)=L(\overline{g_p})
+\alpha_{\overline{g_p}}\c(g_{p+1},\gamma_{p+1})=\c(g,\gamma).\end{equation}
But (\ref{cgnl}) and (\ref{cgn}) give $L(g^n)=L(g)$ for $n\geq 1$
--- this is absurd since in particular $L(g)+\alpha_gL(g)=L(g\cdot
g)=L(g^2)=L(g)$ implies $L(g)=0$, which is impossible.
\\

Since conjugates of hyperbolic automorphisms are hyperbolic, parts
(1) and (2) are now clear.

Now suppose that $L_i$ is a regular length function for each $i$,
and that $g,h\in G$. Then with the notation of Claim 35,
if $g_{p+1}$ and $h_{p+1}$ belong to distinct free factors, we have
$\c(g,h)=L(\overline{g_p})$, $g=\overline{g_{p}}\cdot(g_{p+1}\cdots
g_n)$ and $h=\overline{g_{p}}\cdot(h_{p+1}\cdots h_m)$, and thus
$L(g)=L(\overline{g_p})+\alpha_{\overline{g_p}}L(g_{p+1}\cdots g_n)$
and $L(h)=L(\overline{g_p})+\alpha_{\overline{g_p}}L(h_{p+1}\cdots
h_m)$. Otherwise, let $u_{p+1}$ be the element of the same free
factor as $g_{p+1}$ and $h_{p+1}$  satisfying
$g_{p+1}=u_{p+1}g_{p+1}'$, $h_{p+1}=u_{p+1}h_{p+1}'$,
$L(g_{p+1})=L(u_{p+1})+\alpha_{u_{p+1}}L(g_{p+1}')$,
$L(h_{p+1})=L(u_{p+1})+\alpha_{u_{p+1}}L(h_{p+1}')$ and
$L(u_{p+1})=\c(g_{p+1},h_{p+1})$. It is routine to verify that
$\c(g,h)=L(u)$ where $u=\overline{g_p}\cdot u_{p+1}$, and that
$L(g)=L(u)+\alpha_u L(g')$, $L(h)=L(u)+\alpha_u L(h')$, where
$g=ug'$ and $h=uh'$. This shows that $L$ is regular.
\qed\\

The argument given shows that if the actions of $G_i$ are both free
and rigid then so is that of $G$. In fact if the actions of $G_i$
are merely assumed to be rigid and if for each $i$ there exists
$x_i\in X$ which is not fixed by any $g_i\in G_i$, then taking these
points as the basepoints for the length functions $L_i$, it follows
that all cyclically reduced $g\in G$ of length at least 2 (as a word
in the free product) are hyperbolic and rigid. Thus $G$ acts
rigidly.

\begin{corollary}
\be\item Let $G_i$ be an $\ATF(\R)$ group for $i\in I$. Then
$\ast_{i\in I}G_i$ is an $\ATF(\R)$ group.
\item
The class of $\ATF$ groups is closed under free products. \ee

\end{corollary}

\pf Part (1) is immediate from the theorem. For part (2) it
suffices to show that given ordered abelian groups $\Lambda_i$
($i\in I$) there exist $\Lambda$ and embeddings
$h_i:\Lambda_i\to\Lambda$  such that every automorphism of
$\Lambda_i$ extends to one of $\Lambda$, and such that if free
$\alpha^{(i)}$ actions of $G_i$ on $\Lambda_i$-trees are given, the
induced actions of $G_i$ on $\Lambda\otimes_{\Lambda_i}X_i$ are also
free.

Endow $I$ with a linear order. Let $\Lambda$ be the subgroup of
$\prod_{i\in I} \Lambda_i$ consisting of those $(\lambda_i)$ with
well-ordered support. This makes $\Lambda$ an ordered abelian group.
Moreover each $\Lambda_i$ embeds in $\Lambda$ via $\lambda_i\mapsto
(\delta_{ij}\lambda_i)_{j\in I}$ (where $\delta_{ij}$ is the
Kronecker delta function) and automorphisms $\alpha_*$ of
$\Lambda_i$ extend to automorphisms $\bar{\alpha}_*$ in an obvious
way. By Theorem~\ref{basechange}, there is an induced action of
$G_i$ on a $\Lambda$-tree $\bar{X}_i$. Moreover, choosing a
basepoint $x_i\in X_i$, and letting $L_i=L_{x_i}$ and
$\bar{L}_i=L_{\phi_i(x_i)}$ (where $\phi_i$ denotes the embedding of
$X_i$ in $\bar{X}_i$), we have $\bar{L}_i=h_i\ L_i$, by
Theorem~\ref{basechange}(3)(c). The ancillary functions with respect
to the actions on $X_i$ and $\bar{X}_i$ are similarly related, and
will be similarly notated.

By Proposition~\ref{radius}, $\overline{\Rad}^r(g)$ is spanned by
terms that are $\Z$-linear combinations of $\bar{L}_i=h_i\ L$ and
$\bar{\a}_{i}=h_i\ \a_{i}$. It follows that $\overline{\Rad}^r(g)$
is spanned by $h_i(\Rad^r(g))$.

If $\bar{\lambda}=(\lambda_j)_{j\in I}\in\Lambda$ satisfies
$(1-\bar{\alpha}_g)(\bar{\lambda})=\bar{\b}_i(g)$ then
$\lambda_i-\alpha_g^{(i)}(\lambda_i)=\b_i(g)$. If further
$\bar{\lambda}\in\overline{\Rad}^r(g)$, then $\lambda_i\leq
L_i(\gamma)-\a_i(g)-\alpha_\gamma\a_i(\gamma^{-1}g\gamma)$ for some
$\gamma$. It follows that $\lambda_i\in\Rad^r(g)$, so that $g$ has a
fixed point by Proposition~\ref{hyper-b}, a contradiction. Therefore
$(1-\bar{\alpha}_g^{(i)})^{-1}(\bar{\b}_i(g))\cap\overline{\Rad}^r(g)=\emptyset$,
so that hyperbolicity is preserved by the base change functor in
this case. By Corollary~\ref{free-criterion}(1), the action of $G_i$
on $\bar{X}_i$ is free.

The result now follows from Theorem~\ref{freeproduct}.\qed

\subsection{Ultraproducts}

\begin{theorem}\label{ultraproducts}
Let $G_i$ be a group with an $\alpha^{(i)}$-affine action on a
$\Lambda_i$-tree $(X_i,d_i)$ and let $\mathcal{D}$ be an ultrafilter
in $I$. Then $G=\prod_{i\in I}G_i/\mathcal{D}$ has an induced
$\alpha$-affine action on a $\Lambda$-tree ($\Lambda=\prod_{i\in
I}\Lambda_i/\mathcal{D}$).

If, for almost all $i$, the given actions of $G_i$ are respectively
\be\item free (and without inversions), \item rigid,
\item or regular\ee then so is that of $G$.
\end{theorem}

\pf Much of the proof is a routine use of ultraproducts, and we will
omit most of the details. (See \cite[\S5.5]{Chiswell-book} or
\cite[\S3]{Bell-Slomson} for the necessary background.)

Define a $\Lambda$-tree $(X,d)$ as follows. The set $X$ is the
ultraproduct of the $X_i$, so $X$ consists of $\langle x_i\rangle$
where $x_i\in X_i$ for all $i$, and $d(\langle x_i\rangle,\langle
y_i\rangle)=\langle d_i(x_i,y_i)\rangle$. (Here we denote by
$\langle x_i\rangle$ the equivalence class of the $I$-sequence
$(x_i)_{i\in I}$ where two $I$-sequences are equivalent if their
agreement set is an element of $\mathcal{D}$; similar notational
conventions apply to the other ultra-objects we consider.)

Using the criterion $$y\in[x,z]\Leftrightarrow
d(x,z)=d(x,y)+d(y,z),$$ it is easy to verify that
\begin{equation}\label{ultra-segment} [\langle x_i\rangle,\langle
z_i\rangle]=\langle [x_i,z_i]\rangle.\end{equation}If $\zeta_i$ is
an isometric isomorphism $[0,\lambda_i]\to[x_i,z_i]$, then
$\zeta:[0,\lambda]\to[\langle x_i\rangle,\langle z_i\rangle]$
defined by $$\zeta\langle t_i\rangle=\langle \zeta_i(t_i)\rangle$$
is an isometric isomorphism (where $\lambda=\langle
\lambda_i\rangle$). Therefore $(X,d)$ is geodesic. It is similarly
straightforward to check that for $v\in X$, we have $(x\cdot
y)_v\in\Lambda$, and that $X$ is 0-hyperbolic.

Putting
$$\alpha_{\langle g_i \rangle}\langle \lambda_i\rangle=\langle
\alpha^{(i)}_{g_i}(\lambda_i)\rangle,$$
we obtain a homomorphism $G\to\Aut^+(\Lambda)$
and there is an induced $\alpha$-affine action of $G$ on $X$ given
by $$\langle g_i\rangle\cdot\langle x_i\rangle=\langle g_i
x_i\rangle.$$

If the given actions of $G_i$ are free, it is easy to see that the
action of $G$ is. Using the description of segments in $X$
in~(\ref{ultra-segment}) one can show that rigidity is likewise
preserved by the ultraproduct action. To see that regularity is
preserved by the action of $G$, choose a basepoint $v_i\in X_i$, put
$L_i(g_i)=d_i(v_i,g_iv_i)$, and suppose that for $g_i,h_i\in G_i$
there exists $u_i\in G_i$ such that $\c(g_i,h_i)=L_i(u_i)$. Suppose
further that $L_i(g_i)=L_i(u_i)+\alpha_{u_i}L_i(g_i')$ and
$L_i(h_i)=L_i(u_i)+\alpha_{u_i}L_i(h_i')$ for some $g_i'$ and $h_i'$
such that $g_i=u_ig_i'$ and $h_i=u_ih_i'$. Then one has $\c(\langle
g_i\rangle,\langle h_i\rangle)=L(\langle u_i\rangle)$, where
$L=L_{\langle v_i\rangle}$,  $L\langle g_i\rangle=L\langle
u_i\rangle+\alpha_{\langle u_i\rangle}L\langle g_i'\rangle$,
$L\langle h_i\rangle=L\langle u_i\rangle+\alpha_{\langle
u_i\rangle}L\langle h_i'\rangle$, and $\langle u_i\rangle\langle
g_i'\rangle=\langle g_i\rangle$ and $\langle u_i\rangle\langle
h_i'\rangle=\langle h_i\rangle$. \qed

\begin{theorem}
\be
\item
A group is locally in $\ATF$ if and only if it is in $\ATF$.
\item
A group is fully residually in $\ATF$ if and only if it is in
$\ATF$. \ee\end{theorem}

\pf The necessity in both parts is clear, so we will focus on the
sufficiency.

(1) Let $G_i$ ($i\in I$) denote the finitely generated subgroups of
$G$, and assume that each $G_i$ is in $\ATF$. By
Theorem~\ref{ultraproducts}, the group $^\ast G=\prod_{i\in
I}G_i/\mathcal{D}$ is in $\ATF$, so it suffices to show that
$\mathcal{D}$ can be chosen such that $G$ embeds in $^\ast G$. Write
$a_x=\{i\in I:x\in G_i\}$ for $x\in G$, and note that since
$i_0\in\bigcap_{i=1}^n a_{x_i}$ where $G_{i_0}=\langle
x_1,\ldots,x_n\rangle$, the sets $a_x$ ($x\in G$) have the finite
intersection property. Using Zorn's Lemma, we can find a maximal
ultrafilter $\mathcal{D}$ that contains each $a_x$. Now given $x\in
G$, and putting $x_i=x$ for all $i$ for which $x_i\in G_i$ (and
$x_i=1$ otherwise), we see that $x\mapsto(x_i)_{i\in I}/\mathcal{D}$
is a well-defined monomorphism.

(2) Consider now the set of normal subgroups $N_i$ ($i\in I$) of
$G$  for which $G/N_i$ is in $\ATF$, and suppose that for all finite
subsets $X$ of $G\backslash 1$ we have $N_i\cap X=\emptyset$ for
some $i$. Again we propose to embed $G$ in a suitable ultraproduct
of $\ATF$ groups, from which the required result can be deduced from
Theorem~\ref{ultraproducts}

For $x\in G\backslash 1$ put $a_x=\{i\in I:x\notin N_i\}$. As in
part (1), the sets $a_x$ are easily seen to have the finite
intersection property, whence an ultrafilter $\mathcal{D}$
containing each $a_x$. For $x\in G$ we put $x_i=x$ for all $i$ and
map $x\mapsto(x_iN_i)_{i\in I}/\mathcal{D}$; this gives an embedding
of $G$ in $\prod_{i\in I}(G/N_i)/\mathcal{D}$. \qed

\subsection{Extending an isometric subaction to an affine action}

The following theorem and proof follow the idea of \cite[Corollaire
2(2)]{Liousse}

Recall that an isometric action of $G$ is \emph{abelian} if
$\ell(gh)\leq\ell(g)+\ell(h)$ for all $g,h\in G$.

\begin{theorem}\label{extn-isom}
Let $G$ be a group, and $N$ a normal subgroup. Suppose that $N$ has
a minimal non-abelian isometric action on a $\Lambda$-tree $X$. Let
$\alpha:G\to\Aut^+(\Lambda)$ be a homomorphism, with
$N\leq\ker\alpha$. Then there is an $\alpha$-affine action of $G$ on
$X$ extending the original action of $N$
 if and only if, for $g\in N$ and $\gamma\in G$, we have
\begin{equation}\label{hyper-length-extn}\ell(\gamma g\gamma^{-1})=\alpha_\gamma\ell(g).\end{equation}
\end{theorem}

\pf The sufficiency is straightforward to show.

Conversely, let $\ell$ denote the hyperbolic length function arising
from the isometric action of $N$. Fix $\gamma\in G$ and consider the
$\Lambda$-tree $(X,d_1)$ where $d_1=\alpha_{\gamma^{-1}}\circ d$ and
the ($d_1$-isometric) action of $N$ given by $g\cdot x=\gamma
g\gamma^{-1}x$. Then $d_1(x,g\cdot x)=\alpha_{\gamma^{-1}}d(x,\gamma
g\gamma^{-1}x)$. It follows, using (\ref{hyper-length-extn})
denoting the associated hyperbolic length function by $\ell_1$, that
$\ell_1(g)=\ell(g)$, whence a unique $N$-equivariant isometry
$\phi=\phi_\gamma:(X,d)\to(X,d_1)$, by \cite[Theorem
3.4.1]{Chiswell-book}. Thus $d(\phi(x),\phi(y))=\alpha_\gamma
d(x,y)$, which means that $\phi$ is an affine automorphism of
$(X,d)$ with dilation factor $\alpha_\gamma$.

Now consider $\gamma_1,\gamma_2\in G$. Since $\phi_{\gamma_1}$ is an
$N$-equivariant isometry $(X,d)\to(X,\alpha_{\gamma_1^{-1}}d)$, it
is also an $N$-equivariant isometry
$(X,\alpha_{\gamma_2^{-1}}d)\to(X,\alpha_{\gamma_2^{-1}}\alpha_{\gamma_1^{-1}}d)$.
Therefore $\phi_{\gamma_1}\phi_{\gamma_2}$ is an $N$-equivariant
isometry $(X,d)\to(X,\alpha_{\gamma_2^{-1}\gamma_1^{-1}})$. Since
such an isometry is unique, we must have
$\phi_{\gamma_1\gamma_2}=\phi_{\gamma_1}\phi_{\gamma_2}$.

We therefore have an affine action of $G$ on $X$ by putting $\gamma
x=\phi_\gamma(x)$. Moreover if $\gamma\in N$, then $x\mapsto \gamma
x$ is an $N$-equivariant isometry $(X,d)\to(X,d)$ (since
$\alpha_\gamma$ is trivial in this case), whence
$\phi_\gamma=\gamma$. This means that the action of $G$ extends the
action of $N$.\qed\\

We remark that Chiswell has recently shown that $\ITF(\R^n)$ groups
are right orderable \cite{Chiswell-ro}. Since $\Aut^+(\Lambda)$ is
right orderable, it follows (via the map $\alpha$) that
$\ATF(\Lambda)$ groups are $\ITF(\Lambda)$-by-(right orderable);
hence $\ATF(\R^n)$ groups are right orderable.

\bibliographystyle{plain}
\def\cprime{$'$}

\vspace{0.5cm}

\begin{minipage}[t]{3 in}
\noindent Shane O Rourke\\ Department of Mathematics\\
Cork Institute of Technology\\
Rossa Avenue\\ Cork\\ IRELAND
\\ \verb"shane.orourke@cit.ie"
\end{minipage}

\end{document}